\documentclass{art_ar}

\usepackage{amssymb}

\usepackage{pstricks, pst-plot, pst-node, pst-tree}
\usepackage{array}
\usepackage{fancyhdr}
\usepackage{longtable}
\usepackage{tabularx}
\usepackage{rotating}

\usepackage{makeidx}
\usepackage{amsmath}
\usepackage{theorem}
\usepackage[mathscr]{eucal}
\usepackage{enumerate}
\usepackage[dvips]{epsfig}
\usepackage{float}  % für die Option [H] (genau hier) bei picture usw.

\theoremstyle{plain}
\newtheorem{Def}{Definition}[section]
\newtheorem{Sat}[Def]{Proposition}
\newtheorem{The}[Def]{Theorem}
\newtheorem{Kor}[Def]{Corollary}
\newtheorem{Bem}[Def]{Remark}
\newtheorem{Lem}[Def]{Lemma}

\newcommand{\Id}{\operatorname{Id}}

\begin{document}

\begin{frontmatter}

\title{Stochastic Taylor expansions for the expectation of functionals of diffusion processes}
\author{Andreas R{\"o}{\ss}ler}
\ead{roessler@mathematik.tu-darmstadt.de}
\address{Darmstadt University of Technology, Fachbereich Mathematik, Schlossgartenstr.7,
D-64289 Darmstadt, Germany}

\begin{abstract}
Stochastic Taylor expansions of the expectation of functionals
applied to diffusion processes which are solutions of stochastic
differential equation systems are introduced. Taylor formulas
w.r.t.\ increments of the time are presented for both, It{\^o} and
Stratonovich stochastic differential equation systems with
multi-dimensional Wiener processes. Due to the very complex
formulas arising for higher order expansions, an advantageous
graphical representation by coloured trees is developed. The
convergence of truncated formulas is analyzed and estimates for
the truncation error are calculated. Finally, the stochastic
Taylor formulas based on coloured trees turn out to be a
generalization of the deterministic Taylor formulas using plain
trees as recommended by Butcher for the solutions of ordinary
differential equations.
\end{abstract}

\end{frontmatter}

\section{Introduction} \label{Introduction}
Taylor expansions are well known in the deterministic setting.
They play an important role, e.g. for the determination of
numerical schemes for the approximation of solutions of
differential equations. Especially the modern theory of order
conditions for Runge-Kutta methods due to Butcher~\cite{Butcher87}
is based upon expansions w.r.t.\ trees. A generalization of the
deterministic Taylor formula for a class of It{\^o} processes has
been introduced by Platen and Wagner~\cite{PlWa82}. They make use
of hierarchical sets in order to specify the Taylor expansion of a
functional $f$ of the solution process $X=(X_t)_{t \geq t_0}$ of a
stochastic differential equation (SDE) (see also~\cite{KP99}). In
contrast to the expansion of $f(X_t)$, a Taylor expansion of the
expectation $E(f(X_t))$ is presented in the following. Instead of
using hierarchical sets, coloured rooted trees are used for the
calculation of the Taylor expansion. For example, a representation
with the aid of coloured trees is the key to the stochastic
numerical analysis, especially to the determination of Runge-Kutta
methods for SDEs as investigated for strong approximation by
Burrage and Burrage~\cite{BuBu00a,BuBu96} and for weak
approximation by R{\"o}{\ss}ler~\cite{Roe03}. \\ \\
Let $(\Omega, \mathcal{A}, P)$ be a probability space with a
filtration $(\mathcal{F}_t)_{t \geq 0}$ and let $I=[t_0,T]$ for
some $0 \leq t_0 < T < \infty$. We consider the solution $(X_t)_{t
\in I}$ of a $d$-dimensional autonomous stochastic differential
equation
\begin{equation} \label{Intro-Ito-St-SDE1-autonom-Wm}
    d X_t = a(X_s) \, ds + b(X_s) \, * dW_s, \qquad X_{t_0} = x_0,
\end{equation}
where $W=((W_t^1, \ldots, W_t^m))_{t \geq 0}$ is an
$m$-dimensional Wiener process w.r.t.\ $(\mathcal{F}_t)_{t \geq
0}$ and the constant $x_0 \in \mathbb{R}^d$ is the initial value.
Like in the deterministic setting, each SDE system can be written
as an autonomous SDE system with one additional equation
representing time. Hence without loss of generality, autonomous
stochastic differential equations are treated in the following
only. SDE~(\ref{Intro-Ito-St-SDE1-autonom-Wm}) can be written in
integral form as
\begin{equation} \label{Intro-Ito-St-SDE1-integralform-Wm}
    X_t = x_0 + \int_{t_0}^t a(X_s) \, ds + \sum_{j=1}^m
    \int_{t_0}^t b^j(X_s) \, * dW_s^j
\end{equation}
for $t \in I$, where the $j$th column of the $d \times m$-matrix
function $b=(b^{i,j})$ is denoted by $b^j$ for $j=1, \ldots,m$.
The second integral w.r.t.\ the Wiener process can be interpreted
either as an It{\^o} or a Stratonovich integral which is indicated
by the asterisk. Thus $* dW$ stands for $dW$ in case of It{\^o}
calculus and for $\circ dW$ in case of Stratonovich calculus. We
suppose that $a : \mathbb{R}^d \rightarrow \mathbb{R}^d$ and $b :
\mathbb{R}^d \rightarrow \mathbb{R}^{d \times m}$ are
Borel-measurable and satisfy the conditions of the Existence and
Uniqueness Theorem
(see, e.g., \cite{Arn73,IW89,KS99,Ok98}). \\ \\
In order to investigate the Taylor expansion of the expectation of
some functional $f : \mathbb{R}^d \rightarrow \mathbb{R}$ of the
solution $X$ of SDE~(\ref{Intro-Ito-St-SDE1-integralform-Wm}) in
the case of It{\^o} calculus, we will make use of the diffusion
operator $L^0$ and the operators $L^j$ for $j=1, \ldots, m$ in
case of an $m$-dimensional Wiener process. They are defined as
\begin{equation} \label{Operator_L0_L1}
    \begin{split}
    L^0 = &\sum_{k=1}^d a^k \, \frac{\partial}{\partial x^k} +
    \frac{1}{2} \sum_{k,l=1}^d \sum_{j=1}^m b^{k,j} \,  b^{l,j} \,
    \frac{\partial^2}{\partial x^k \partial x^l}
    \\
    L^j = &\sum_{k=1}^d b^{k,j} \, \frac{\partial}{\partial x^k}
    \qquad j=1, \ldots, m.
    \end{split}
\end{equation}
Applying It{\^o}'s formula recursively to $f(X_t)$ and afterwards
to $L^0 f(X_s)$ yields
\begin{equation} \label{Taylor-recursive-intro:1}
    \begin{split}
    f(X_t) = &f(X_{t_0}) + \int_{t_0}^t L^0 f(X_s) \, ds +
    \sum_{j=1}^m \int_{t_0}^t L^j f(X_s) \, dW_s^j \\
    = &f(X_{t_0}) + \int_{t_0}^t \big( L^0 f(X_{t_0}) +
    \int_{t_0}^s L^0 L^0 f(X_u) \, du \\
    &+ \sum_{j=1}^m \int_{t_0}^s L^j L^0 f(X_u) \, dW_u^j
    \big) \, ds + \sum_{j=1}^m \int_{t_0}^t L^j f(X_s) \, dW_s^j
    \end{split}
\end{equation}
This procedure can be continued by applying again It{\^o}'s
formula to $L^0 L^0 f(X_u)$ and so on. In the following, it is
always assumed that all of the necessary derivatives and multiple
integrals exist.
If we consider the case of $d=m=1$ and take the expectation, then
we get with $X_{t_0} = x_0 \in \mathbb{R}^d$
\begin{equation}
\begin{split} \label{Intro-Ito-Taylor-d1}
    E^{t_0,x_0}(f(X_t)) = & f(x_0) + \left( f' \, a +
    \tfrac{1}{2} f'' \, b^2 \right)(x_0) \cdot (t-t_0)
    + \left( f' \cdot \left( a \, a' +
    \tfrac{1}{2} b^2 \, a'' \right) \right. \\
    &+ f'' \cdot \left( a^2 + a \, b \, b' + b^2 \, a'
    + \tfrac{1}{2} b^2 \, {b'}^2 +
    \tfrac{1}{2} b^3 \, b'' \right) \\
    & \left. + f''' \cdot \left( a \, b^2 + b^3 \, b' \right)
    + f'''' \cdot \left( \tfrac{1}{4}
    b^4 \right) \right)(x_0) \cdot \frac{(t-t_0)^2}{2!} + \ldots
\end{split}
\end{equation}
since the expectation of a multiple It{\^o} integral with at least
one integration w.r.t.\ a Wiener process
vanishes~\cite{IW89,KP99}. The terms on the right hand side, which
are called elementary differentials in numerical analysis, get
more and more complex for expansions of higher orders. Thus, an
approach similar to the rooted tree theory for deterministic
ordinary differential equations as described by
Butcher~\cite{Butcher87} is helpful.
\section{Rooted Tree Theory} \label{Rooted-Tree-Theory}
\begin{figure}[H]
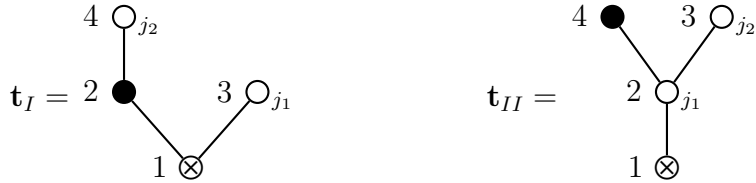

\begin{center}
\begin{tabular}{ccc}
    $\textbf{t}_I = \begin{array}{c}
    \text{\pstree[treemode=U, dotstyle=otimes, dotsize=3.2mm, levelsep=0.1cm, radius=1.6mm, treefit=loose]
    {\Tn}{
    \pstree[treemode=U, dotstyle=otimes, dotsize=3.2mm, levelsep=1cm, radius=1.6mm, treefit=loose]
    {\Tdot~[tnpos=l]{1 }}
    {\pstree{\TC*~[tnpos=l]{2}}{\TC~[tnpos=l]{4}~[tnpos=r]{$\!\!_{j_2}$}}
    \TC~[tnpos=l]{3}~[tnpos=r]{$\!\!_{j_1}$}}
    }}
    \end{array}$
    & \qquad \qquad \qquad &
    $\textbf{t}_{II} = \begin{array}{c}
    \text{\pstree[treemode=U, dotstyle=otimes, dotsize=3.2mm, levelsep=0.1cm, radius=1.6mm, treefit=loose]
    {\Tn}{
    \pstree[treemode=U, dotstyle=otimes, dotsize=3.2mm, levelsep=1cm, radius=1.6mm, treefit=loose]
    {\Tdot~[tnpos=l]{1 }} {\pstree{\TC~[tnpos=l]{2}~[tnpos=r]{$\!\!_{j_1}$}}{\TC*~[tnpos=l]{4}
    \TC~[tnpos=l]{3}~[tnpos=r]{$\!\!_{j_2}$}}}
    }}
    \end{array}$
\end{tabular}
\caption{Two elements of $LTS$ with $j_1, j_2 \in \{1, \ldots,
m\}$.} \label{St-S-tree-examples-tI+tII:Wm}
\end{center}
\end{figure}

\begin{figure}[htb]
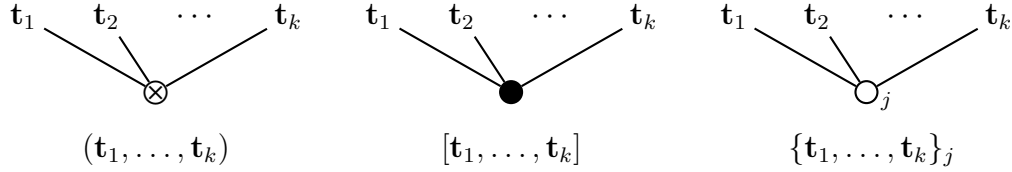

\begin{center}
\begin{tabular}{ccccc}
    \pstree[treemode=U, dotstyle=otimes, dotsize=3.2mm, levelsep=0.1cm, radius=1.6mm, treefit=loose]
    {\Tn}{
    \pstree[treemode=U, dotstyle=otimes, dotsize=3.2mm, levelsep=1cm, radius=1.6mm, treefit=loose, nodesepB=1mm]
    {\Tdot} {\Tr{$\textbf{t}_1$} \Tr{$\textbf{t}_2$} \Tr[edge=none]{$\cdots$} \Tr{$\textbf{t}_k$}}
    }
    & &
    \pstree[treemode=U, dotstyle=otimes, dotsize=3.2mm, levelsep=0.1cm, radius=1.6mm, treefit=loose]
    {\Tn}{
    \pstree[treemode=U, dotstyle=otimes, dotsize=3.2mm, levelsep=1cm, radius=1.6mm, treefit=loose, nodesepB=1mm]
    {\TC*} {\Tr{$\textbf{t}_1$} \Tr{$\textbf{t}_2$} \Tr[edge=none]{$\cdots$} \Tr{$\textbf{t}_k$}}
    }
    & &
    \pstree[treemode=U, dotstyle=otimes, dotsize=3.2mm, levelsep=0.1cm, radius=1.6mm, treefit=loose]
    {\Tn}{
    \pstree[treemode=U, dotstyle=otimes, dotsize=3.2mm, levelsep=1cm, radius=1.6mm, treefit=loose, nodesepB=1mm]
    {\TC~[tnpos=r]{$\!\!_{j}$}} {\Tr{$\textbf{t}_1$} \Tr{$\textbf{t}_2$} \Tr[edge=none]{$\cdots$} \Tr{$\textbf{t}_k$}}
    }
    \\
    $(\textbf{t}_1, \ldots, \textbf{t}_k)$ & &
    $[\textbf{t}_1, \ldots, \textbf{t}_k]$ & & $\,\, \{\textbf{t}_1, \ldots, \textbf{t}_k\}_j$
\end{tabular}
\caption{Writing a coloured S-tree with brackets.}
\label{St-tree-bracket-together}
\end{center}
\end{figure}

\begin{figure}[H]
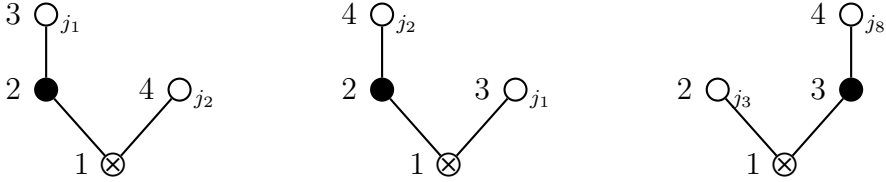

\begin{center}
    \begin{tabular}{ccccc}
    \pstree[treemode=U, dotstyle=otimes, dotsize=3.2mm, levelsep=0.1cm, radius=1.6mm, treefit=loose]
    {\Tn}{
    \pstree[treemode=U, dotstyle=otimes, dotsize=3.2mm, levelsep=1cm, radius=1.6mm, treefit=loose]
    {\Tdot~[tnpos=l]{1 }} {\pstree{\TC*~[tnpos=l]{2}}{\TC~[tnpos=l]{3}~[tnpos=r]{$\!\!_{j_1}$}}
    \TC~[tnpos=l]{4}~[tnpos=r]{$\!\!_{j_2}$}}
    }
    & \qquad \qquad &
    \pstree[treemode=U, dotstyle=otimes, dotsize=3.2mm, levelsep=0.1cm, radius=1.6mm, treefit=loose]
    {\Tn}{
    \pstree[treemode=U, dotstyle=otimes, dotsize=3.2mm, levelsep=1cm, radius=1.6mm, treefit=loose]
    {\Tdot~[tnpos=l]{1 }} {\pstree{\TC*~[tnpos=l]{2}}{\TC~[tnpos=l]{4}~[tnpos=r]{$\!\!_{j_2}$}}
    \TC~[tnpos=l]{3}~[tnpos=r]{$\!\!_{j_1}$}}
    }
    & \qquad \qquad &
    \pstree[treemode=U, dotstyle=otimes, dotsize=3.2mm, levelsep=0.1cm, radius=1.6mm, treefit=loose]
    {\Tn}{
    \pstree[treemode=U, dotstyle=otimes, dotsize=3.2mm, levelsep=1cm, radius=1.6mm, treefit=loose]
    {\Tdot~[tnpos=l]{1 }} {\TC~[tnpos=l]{2}~[tnpos=r]{$\!\!_{j_3}$} \pstree{\TC*~[tnpos=l]{3}}{
    \TC~[tnpos=l]{4}~[tnpos=r]{$\!\!_{j_8}$}}}
    }
    \end{tabular}
\caption{Trees of the same equivalence class.}
\label{St-equal-trees:Wm}
\end{center}
\end{figure}

First of all a definition of the coloured graphs which will be
suitable in the rooted tree theory for stochastic differential
equation systems is given. In contrast to the trees introduced by
Burrage and Burrage~\cite{BuBu00a,BuBu96}, the trees defined in
the following own an additional kind of node corresponding to the
functional and have a particular structure.
\begin{Def} \label{Def:rooted-S-trees:Wm}
    Let $l$ be a positive integer.
    \begin{enumerate}
        \item A {\emph{monotonically labelled S-tree (stochastic
        tree)}} $\textbf{t}$
        with $l=l(\textbf{t})$ nodes is a pair of maps $\textbf{t}=(\textbf{t}',\textbf{t}'')$
        with $\textbf{t}''=(\textbf{t}_A'',\textbf{t}_J'')$ and
        \begin{equation*}
            \begin{split}
                \textbf{t}' & : \{2, \ldots, l\} \longrightarrow \{1, \ldots,
                l-1\} \\
                \textbf{t}_A'' & : \{1, \ldots, l\} \longrightarrow A \\
                \textbf{t}_J'' & : \{1, \ldots, l\} \longrightarrow J \cup
                \{ 0 \}
            \end{split}
        \end{equation*}
        so that $\textbf{t}'(i) < i$ for $i=2, \ldots, l$; $A$ denotes
        a finite set and $J=\{j_1, \ldots, j_s\}$ is a set of indices.
        \item $LTS$ denotes the set of all monotonically labelled
        S-trees. Here two trees $\textbf{t}=(\textbf{t}',\textbf{t}'')$ and $\textbf{u}=(\textbf{u}',\textbf{u}'')$
        just differing by their renamed indices $\textbf{t}_J''$ and $\textbf{u}_J''$ are
        considered to be identical if there exists a bijective
        map
        \begin{equation}
        \pi : \{j_1, \ldots, j_s\} \cup \{ 0 \}
        \longrightarrow
        \{j_{1}, \ldots, j_{s}\} \cup \{ 0 \} \quad \text{
        with } \pi(0) = 0
        \end{equation}
        so that $\textbf{t}_J''(i) = \pi(\textbf{u}_J''(i))$ holds for $i=1, \ldots,
        l$.
    \end{enumerate}
\end{Def}
\begin{Bem}
    Thus $LTS$ consists of monotonically labelled trees whereby two
    trees are assumed to be identical if we can rename the indices $j_i$ of the
    first one to the indices of the second tree, so that nodes
    having identical indices receive again identical indices and
    nodes with different indices receive again different indices.
\end{Bem}
Every monotonically labelled S-tree $\textbf{t}$ can be
represented as a graph, whose nodes are elements of $\{1, \ldots,
l(\textbf{t})\}$ and whose arcs are the pairs $(\textbf{t}'(i),i)$
for $i=2, \ldots, l(\textbf{t})$. Here, $\textbf{t}'$ defines a
father son relation between the nodes, i.e.\ $\textbf{t}'(i)$ is
the father of the son $i$. Furthermore the colour
$\textbf{t}''(i)=(\textbf{t}_A''(i), \textbf{t}_J''(i))$, which
consists of one element of the set $A$ and one element of the set
$J \cup \{ 0 \}$, is added to the node $i$ for $i=1, \ldots,
l(\textbf{t})$. The node with label 1 is called {\emph{root}} and
it is always sketched as the lowest node of the graph.
If not stated otherwise, we take $A = \{ \gamma, \tau, \sigma \}$
in the following. In this case $\gamma = \,\,$~\pstree[treemode=U,
dotstyle=otimes, dotsize=3.2mm, levelsep=0.1cm, radius=1.6mm,
treefit=loose]
    {\Tn}{
    \pstree[treemode=U, dotstyle=otimes, dotsize=3.2mm, levelsep=0cm, radius=1.6mm, treefit=loose]
    {\Tdot~[tnpos=r]{ }} {}
    }
denotes the root, $\tau = $~\pstree[treemode=U, dotstyle=otimes,
dotsize=3.2mm, levelsep=0.1cm, radius=1.6mm, treefit=loose]
    {\Tn}{
    \pstree[treemode=U, dotstyle=otimes, dotsize=3.2mm, levelsep=0cm, radius=1.6mm, treefit=loose]
    {\TC*~[tnpos=r]{}} {}
    }
is a deterministic node and $\sigma = $~\pstree[treemode=U,
dotstyle=otimes, dotsize=3.2mm, levelsep=0.1cm, radius=1.6mm,
treefit=loose]
    {\Tn}{
    \pstree[treemode=U, dotstyle=otimes, dotsize=3.2mm, levelsep=0cm, radius=1.6mm, treefit=loose]
    {\TC~[tnpos=r]{}} {}
    }
is a stochastic node. Further we consider a set $J=\{j_1, \ldots,
j_s \}$ of indices which are taken as variables so that $j_i \in
\{1, \ldots, m\}$ for $i=1, \ldots, s$. Usually we assign to the
root $\gamma$ and to all deterministic nodes $\tau$ the index $0$,
which will be omitted in the following. However, we assign to
every stochastic node $\sigma$ an index $j_i$, which is associated
with the $j_i$th component of the corresponding $m$-dimensional
Wiener process of the considered SDE. In case of a one-dimensional
Wiener process one can omit the set of indices $J$ since $j_i = 1$
for all $i=1, \ldots, s$ (see also \cite{Roe03}). As an example
Figure~\ref{St-S-tree-examples-tI+tII:Wm} presents two elements of
$LTS$.
\noindent For the labelled S-tree $\textbf{t}_I$ in
Figure~\ref{St-S-tree-examples-tI+tII:Wm} we have
$\textbf{t}_I'(2)=\textbf{t}_I'(3)=1$ and $\textbf{t}_I'(4)=2$.
The colour of the nodes is given by $\textbf{t}_I''(1)=\gamma$,
$\textbf{t}_I''(2)=\tau$, $\textbf{t}_I''(3)=\sigma_{j_1}$ and
$\textbf{t}_I''(4)=\sigma_{j_2}$.
\begin{Def} \label{Def:order-S-tree-W1}
    Let $\textbf{t}=(\textbf{t}',\textbf{t}'') \in LTS$ w.r.t.\ the set $A=\{\gamma, \tau, \sigma\}$.
    Then we denote by $d(\textbf{t}) = \sharp \{ i : \textbf{t}''(i) = \tau \}$
    the number of deterministic nodes and by
    $s(\textbf{t}) = \sharp \{ i : \textbf{t}''(i) = \sigma_j, \, j \in J \}$ the number
    of stochastic nodes of $\textbf{t}$.
    The {\emph{order}} $\rho(\textbf{t})$ of the tree $\textbf{t}$ is defined as
    $\rho(\textbf{t}) = d(\textbf{t}) + \tfrac{1}{2} s(\textbf{t})$ with $\rho(\gamma) = 0$.
\end{Def}
The order of the trees $\textbf{t}_I$ and $\textbf{t}_{II}$
presented in Figure~\ref{St-S-tree-examples-tI+tII:Wm} can be
calculated as $\rho(\textbf{t}_I)=\rho(\textbf{t}_{II})=2$. Every
labelled S-tree can be written as a combination of three different
brackets defined as follows.
\begin{Def}
    If $\textbf{t}_1, \ldots, \textbf{t}_k$ are coloured trees then we denote by
    \begin{equation*}
        (\textbf{t}_1, \ldots, \textbf{t}_k), \,\,\,\,\,
        [\textbf{t}_1, \ldots, \textbf{t}_k] \,\,\,\,\, \text{ and }
        \,\,\,\,\, \{\textbf{t}_1, \ldots, \textbf{t}_k \}_j
    \end{equation*}
    the tree in which $\textbf{t}_1, \ldots, \textbf{t}_k$ are each joined by a
    single branch to $\,\,$
    \pstree[treemode=U, dotstyle=otimes, dotsize=3.2mm, levelsep=0.1cm, radius=1.6mm, treefit=loose]
    {\Tn}{
    \pstree[treemode=U, dotstyle=otimes, dotsize=3.2mm, levelsep=0cm, radius=1.6mm, treefit=loose]
    {\Tdot} {}
    }
    $\,$,
    \pstree[treemode=U, dotstyle=otimes, dotsize=3.2mm, levelsep=0.1cm, radius=1.6mm, treefit=loose]
    {\Tn}{
    \pstree[treemode=U, dotstyle=otimes, dotsize=3.2mm, levelsep=0cm, radius=1.6mm, treefit=loose]
    {\TC*} {}
    }
    $\,\,$and
    \pstree[treemode=U, dotstyle=otimes, dotsize=3.2mm, levelsep=0.1cm, radius=1.6mm, treefit=loose]
    {\Tn}{
    \pstree[treemode=U, dotstyle=otimes, dotsize=3.2mm, levelsep=0cm, radius=1.6mm, treefit=loose]
    {\TC~[tnpos=r]{$\!\!_j$}} {}
    },
    respectively (see Figure~\ref{St-tree-bracket-together}).
\end{Def}
Therefore proceeding recursively, for the two examples
$\textbf{t}_I$ and $\textbf{t}_{II}$ in
Figure~\ref{St-S-tree-examples-tI+tII:Wm} we obtain
    $\textbf{t}_I = ([
    \text{\pstree[treemode=U, dotstyle=otimes, dotsize=3.2mm, levelsep=0.1cm, radius=1.6mm, treefit=loose]
    {\Tn}{
    \pstree[treemode=U, dotstyle=otimes, dotsize=3.2mm, levelsep=0cm, radius=1.6mm, treefit=loose]
    {\TC~[tnpos=r]{$\!\!_{j_2}$}} {}
    }} ],
    \text{\pstree[treemode=U, dotstyle=otimes, dotsize=3.2mm, levelsep=0.1cm, radius=1.6mm, treefit=loose]
    {\Tn}{
    \pstree[treemode=U, dotstyle=otimes, dotsize=3.2mm, levelsep=0cm, radius=1.6mm, treefit=loose]
    {\TC~[tnpos=r]{$\!\!_{j_1}$}} {}
    }} ) = ([\sigma_{j_2}] , \sigma_{j_1})$
    and
    $\textbf{t}_{II} = (\{
    \text{\pstree[treemode=U, dotstyle=otimes, dotsize=3.2mm, levelsep=0.1cm, radius=1.6mm, treefit=loose]
    {\Tn}{
    \pstree[treemode=U, dotstyle=otimes, dotsize=3.2mm, levelsep=0cm, radius=1.6mm, treefit=loose]
    {\TC*} {}
    }},
    \text{\pstree[treemode=U, dotstyle=otimes, dotsize=3.2mm, levelsep=0.1cm, radius=1.6mm, treefit=loose]
    {\Tn}{
    \pstree[treemode=U, dotstyle=otimes, dotsize=3.2mm, levelsep=0cm, radius=1.6mm, treefit=loose]
    {\TC~[tnpos=r]{$\!\!_{j_2}$}} {}
    }}
    \}_{j_1}) = (\{ \tau, \sigma_{j_2} \}_{j_1})$. \\ \\
%
% ===========================================================
%
In the following, we will concentrate our considerations to one
representative tree of each equivalence class defined as follows:
\begin{Def} \label{St-tree-equivalence:Wm}
    Let $\textbf{t}=(\textbf{t}',\textbf{t}'')$ and $\textbf{u}=(\textbf{u}',\textbf{u}'')$
    be elements of $LTS$. Then
    the trees $\textbf{t}$ and $\textbf{u}$ are equivalent, i.e.\ $\textbf{t} \sim \textbf{u}$, if the
    following hold:
    \begin{enumerate}[(i)]
        \item $l(\textbf{t})=l(\textbf{u})$
        \item There exist two bijective maps
            \begin{equation*}
                \begin{split}
                \psi &: \{1, \ldots, l(\textbf{t})\} \rightarrow \{1,
                \ldots, l(\textbf{t})\} \quad \text{ with } \quad
                \psi(1)=1, \\
                \pi &: \{j_1, \ldots, j_s\} \cup \{ 0 \} \rightarrow
                \{j_{1}, \ldots, j_{s}\} \cup \{ 0 \} \quad \text{
                with } \pi(0) = 0
                \end{split}
            \end{equation*}
            so that the following diagram commutes
            \begin{equation*}
                \begin{psmatrix}[colsep=3cm, rowsep=.5cm]
                    \{2,\ldots,l(\textbf{t})\} & \{1,\ldots,l(\textbf{t})\} & \\
                    & & A \times (J \cup \{ 0 \}) \\
                    \{2,\ldots,l(\textbf{t})\} & \{1,\ldots,l(\textbf{t})\} &
                    \psset{arrows=->, nodesep=3pt}
                    \everypsbox{\scriptstyle}
                    \ncline{1,1}{1,2}\taput{\textbf{t}'}
                    \ncline{3,1}{3,2}\taput{\textbf{u}'}
                    \ncline{1,1}{3,1}\trput{\psi}
                    \ncline{1,2}{3,2}\trput{\psi}
                    \ncline{1,2}{2,3}\taput{(\textbf{t}_A'',\textbf{t}_J'')}
                    \ncline{3,2}{2,3}\tbput{(\textbf{u}_A'',\pi(\textbf{u}_J''))}
                \end{psmatrix}
            \end{equation*}
    \end{enumerate}
\end{Def}
Thus, a monotonically labelled S-tree $\textbf{u}$ is equivalent
to $\textbf{t}$, if each label $k$ is replaced by $\psi(k)$ and if
the variable indices $j_i$ are replaced by $\pi(j_{i})$.
Now it is straight forward to consider equivalence classes.
\begin{Def} \label{Chap2:cardinality:Wm}
    The set of all equivalence classes under the relation of
    Definition~\ref{St-tree-equivalence:Wm} is denoted by $TS = LTS /
    \sim$. The elements of $TS$ are called (rooted)
    {\emph{S-trees}}. We denote by $\alpha(\textbf{t})$ the cardinality of
    $\textbf{t}$, i.e.\ the number of possibilities of monotonically
    labelling the nodes of $\textbf{t}$ with numbers $1, \ldots, l(\textbf{t})$.
\end{Def}
Thus, all trees in Figure~\ref{St-equal-trees:Wm} belong to the
same equivalence class as $\textbf{t}_I$ in the example above,
since the indices $j_1$ and $j_2$ are just renamed either by $j_2$
and $j_1$ or $j_8$ and $j_3$, respectively. Finally the graphs
differ only in the labelling of their number indices.
%
% =======================================================
%
For every rooted tree $\textbf{t} \in LTS$, there exists a
corresponding {\emph{elementary differential}} which is a direct
generalization of the differential in the deterministic case (see,
e.g., \cite{Butcher87}). For $j \in \{1, \ldots, m\}$, the
elementary differential is defined recursively by
\begin{equation*}
    F(\gamma)(x) = f(x), \qquad
    F(\tau)(x) = a(x), \qquad
    F(\sigma_j)(x) = b^j(x),
\end{equation*}
for single nodes and by
\begin{equation} \label{St-elementary-differential-F:Wm}
    F(\textbf{t})(x) =
    \begin{cases}
    f^{(k)}(x) \cdot (F(\textbf{t}_1)(x), \ldots, F(\textbf{t}_k)(x)) &
    \text{for } \textbf{t}=(\textbf{t}_1, \ldots, \textbf{t}_k) \\
    a^{(k)}(x) \cdot (F(\textbf{t}_1)(x), \ldots,
    F(\textbf{t}_k)(x)) & \text{for } \textbf{t}=[\textbf{t}_1, \ldots, \textbf{t}_k] \\
    {b^j}^{(k)}(x) \cdot (F(\textbf{t}_1)(x),
    \ldots, F(\textbf{t}_k)(x)) & \text{for } \textbf{t}=\{\textbf{t}_1, \ldots,
    \textbf{t}_k\}_j
    \end{cases}
\end{equation}
for a tree $\textbf{t}$ with more than one node. Here $f^{(k)}$,
$a^{(k)}$ and ${b^j}^{(k)}$ define a symmetric $k$-linear
differential operator, and one can choose the sequence of labelled
S-trees $\textbf{t}_1, \ldots, \textbf{t}_k$ in an arbitrary
order. For example, the $I$th component of $a^{(k)} \cdot
(F(\textbf{t}_1), \ldots, F(\textbf{t}_k))$ can be written as
\begin{align*}
    ( a^{(k)} \cdot (F(\textbf{t}_1), \ldots, F(\textbf{t}_k)) )^I
    &= \sum_{J_1, \ldots, J_k=1}^d \frac{\partial^k
    a^I}{\partial x^{J_1} \ldots \partial x^{J_k}} \,
    (F^{J_1}(\textbf{t}_1), \ldots, F^{J_k}(\textbf{t}_k))
\end{align*}
where the components of vectors are denoted by superscript
indices, which are chosen as capitals.
As a result of this we get for $\textbf{t}_I$ and
$\textbf{t}_{II}$ the elementary differentials
\begin{equation*}
    \begin{split}
    &F(\textbf{t}_I) = f'' (a'(b^{j_2}), b^{j_1}) = \sum_{J_1,J_2=1}^d
    \frac{\partial^2 f}{\partial x^{J_1} \partial x^{J_2}}
    \big( \sum_{K_1=1}^d \frac{\partial a^{J_1}}{\partial x^{K_1}}
    \, b^{K_1,j_2} \cdot b^{J_2,j_1} \big) \\
    &F(\textbf{t}_{II}) = f' ({b^{j_1}}'' (a, b^{j_2})) = \sum_{J_1=1}^d
    \frac{\partial f}{\partial x^{J_1}} \big( \sum_{K_1, K_2 =1}^d
    \frac{\partial^2 b^{J_1,j_1}}{\partial x^{K_1} \partial
    x^{K_2}} \, a^{K_1} \cdot b^{K_2,j_2} \big)
    \end{split}
\end{equation*}
It has to be pointed out that the elementary differentials for the
trees presented in Figure~\ref{St-equal-trees:Wm} coincide with
$F(\textbf{t}_I)$ if the variable indices $j_i$ are simply renamed
by a suitable bijective mapping $\pi$. \\ \\
There exists a close relation between differentiation of an
elementary differential and the growth of the corresponding
labelled tree, that is adding further nodes to the tree.
\begin{Lem} \label{Lem-tree-growth}
    Let $\textbf{t}=(\textbf{t}',\textbf{t}'') \in LTS$ and let $\lambda = l(\textbf{t})$
    denote the number of nodes of $\textbf{t}$. Then for $g \equiv a$ or $g \equiv
    b^j$ and $x \in \mathbb{R}^d$
    \begin{equation}
        \sum_{k=1}^d g^k(x) \, \frac{\partial}{\partial x^k} \,
        F(\textbf{t})(x) = \sum_{\textbf{u} \in H_1(\textbf{t})} F(\textbf{u})(x)
    \end{equation}
    holds, where $H_1(\textbf{t})$ is the set of trees $\textbf{u}=(\textbf{u}',\textbf{u}'') \in LTS$ with
    $l(\textbf{u})=\lambda+1$ nodes, $\textbf{u}'|_{\{2,\ldots,\lambda\}}=
    \textbf{t}'$, $\textbf{u}''|_{\{1, \ldots, \lambda\}}
    = \textbf{t}''$ and either $\textbf{u}''(\lambda+1)=\tau$ in case of $g \equiv a$
    or $\textbf{u}''(\lambda+1)=\sigma_j$ in case of $g \equiv b^j$.
\end{Lem}
{\bf{Proof.}} $F(\textbf{t})(x)$ is by definition a product of
some derivatives of the functions $f$, $a$ and $b^j$. So if we
apply the operator $\sum_{k=1}^d g^k \frac{\partial}{\partial
x^k}$ to $F(\textbf{t})$, we have to use the product rule for
differentiation. Thus we get a sum with exactly $\lambda =
l(\textbf{t})$ summands. Now, we only have to observe that each of
these summands is equal to the elementary differential of the tree
$\textbf{t}$ with one additional arc $(\lambda+1,i)$ and a new
node $\lambda+1$ of type $\tau$ in case of $g \equiv a$ or of type
$\sigma_j$ in case of $g \equiv b^j$ for the corresponding $i \in
\{1, \ldots, \lambda\}$,
respectively. \hfill $\Box$ \\ \\
The growth of an $S$-tree which corresponds to the twice repeated
differentiation of the elementary differential, is stated by the
following lemma.
\begin{Lem} \label{Lem-tree-growth-2}
    Let $\textbf{t}=(\textbf{t}',\textbf{t}'') \in LTS$ and let $\lambda = l(\textbf{t})$
    denote the number of nodes of $\textbf{t}$. Then for $x \in
    \mathbb{R}^d$
    \begin{equation}
        \sum_{k,l=1}^d b^{k,j}(x) \, b^{l,j}(x) \, \frac{\partial^2}{\partial x^k \partial
        x^l} \, F(\textbf{t})(x) = \sum_{\textbf{u} \in H_2(\textbf{t})} F(\textbf{u})(x)
    \end{equation}
    holds, where $H_2(\textbf{t})$ is the set of trees $\textbf{u}=(\textbf{u}',\textbf{u}'') \in LTS$ with
    $l(\textbf{u})=\lambda+2$ nodes, $\textbf{u}'|_{\{2,\ldots,\lambda\}}=
    \textbf{t}'$, $\textbf{u}''|_{\{1, \ldots, \lambda\}}
    = \textbf{t}''$, $\textbf{u}''(\lambda+1)= \textbf{u}''(\lambda+2)
    = \sigma_j$ and $\textbf{u}'(\lambda+2) \neq \lambda+1$.
\end{Lem}
{\bf{Proof.}} Again, $F(\textbf{t})(x)$ is by definition a product
of some derivatives of the functions $f$, $a$ and $b^j$. So if we
apply the operator $\sum_{k,l=1}^d b^{k,j} \, b^{l,j}
\frac{\partial^2}{\partial x^k \partial x^l}$ to $F(\textbf{t})$,
we have to use the product rule for differentiation. Thus we
obtain a sum with exactly $\lambda^2 = l(\textbf{t})^2$ summands.
Now, we only have to observe that each of these summands is equal
to the elementary differential of the tree $\textbf{t}$ with two
additional arcs $(\lambda+1,r)$ and $(\lambda+2,s)$ and two new
nodes $\lambda+1$ and $\lambda+2$, both of type $\sigma_j$, for
the corresponding $r,s \in \{1, \ldots, \lambda\}$, respectively.
It has to be pointed out, that $r=s$ or $r \neq s$
is possible.\hfill $\Box$ \\ \\
In the following, let $C_P^l(\mathbb{R}^d, \mathbb{R})$ denote the
space of $l$ times continuously differentiable functions $g \in
C^l(\mathbb{R}^d, \mathbb{R})$ for which all partial derivatives
up to order $l$ have {\emph{polynomial growth}}. That is, for
which there exist constants $K > 0$ and $r \in \mathbb{N}$
depending on $g$, such that
\begin{equation} \label{poly_growth}
    | \partial_x^i g(x) | \leq K \, (1+\|x\|^{2r})
\end{equation}
holds for all $x \in \mathbb{R}^d$ and any partial derivative
$\partial_x^i g$ of order $i \leq l$. %\\ \\
\section{Truncated Taylor expansions for It{\^o} SDE Systems}
\begin{figure}[H]
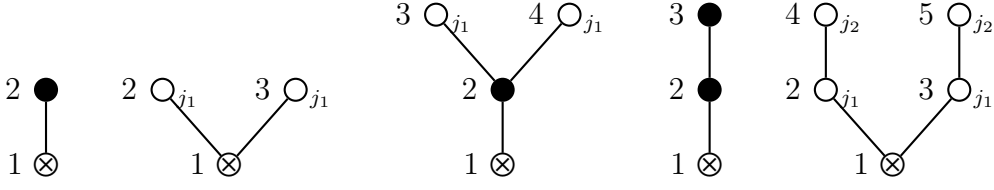

\begin{center}
    \begin{tabular}{ccccccccc}
    \pstree[treemode=U, dotstyle=otimes, dotsize=3.2mm, levelsep=0.1cm, radius=1.6mm, treefit=loose]
    {\Tn}{
    \pstree[treemode=U, dotstyle=otimes, dotsize=3.2mm, levelsep=1cm, radius=1.6mm, treefit=loose]
    {\Tdot~[tnpos=l]{1 }} {\TC*~[tnpos=l]{2}}
    }
    &
    \quad
    &
    \pstree[treemode=U, dotstyle=otimes, dotsize=3.2mm, levelsep=0.1cm, radius=1.6mm, treefit=loose]
    {\Tn}{
    \pstree[treemode=U, dotstyle=otimes, dotsize=3.2mm, levelsep=1cm, radius=1.6mm, treefit=loose]
    {\Tdot~[tnpos=l]{1 }} {\TC~[tnpos=l]{2}~[tnpos=r]{$\!\!_{j_1}$} \TC~[tnpos=l]{3}~[tnpos=r]{$\!\!_{j_1}$}}
    }
    &
    \quad
    &
    \pstree[treemode=U, dotstyle=otimes, dotsize=3.2mm, levelsep=0.1cm, radius=1.6mm, treefit=loose]
    {\Tn}{
    \pstree[treemode=U, dotstyle=otimes, dotsize=3.2mm, levelsep=1cm, radius=1.6mm, treefit=loose]
    {\Tdot~[tnpos=l]{1 }} {\pstree{\TC*~[tnpos=l]{2}}{\TC~[tnpos=l]{3}~[tnpos=r]{$\!\!_{j_1}$}
    \TC~[tnpos=l]{4}~[tnpos=r]{$\!\!_{j_1}$}}}
    }
    &
    \quad
    &
    \pstree[treemode=U, dotstyle=otimes, dotsize=3.2mm, levelsep=0.1cm, radius=1.6mm, treefit=loose]
    {\Tn}{
    \pstree[treemode=U, dotstyle=otimes, dotsize=3.2mm, levelsep=1cm, radius=1.6mm, treefit=loose]
    {\Tdot~[tnpos=l]{1 }} {\pstree{\TC*~[tnpos=l]{2}}{\TC*~[tnpos=l]{3}}}
    }
    &
    \quad
    &
    \pstree[treemode=U, dotstyle=otimes, dotsize=3.2mm, levelsep=0.1cm, radius=1.6mm, treefit=loose]
    {\Tn}{
    \pstree[treemode=U, dotstyle=otimes, dotsize=3.2mm, levelsep=1cm, radius=1.6mm, treefit=loose]
    {\Tdot~[tnpos=l]{1 }} {\pstree{\TC~[tnpos=l]{2}~[tnpos=r]{$\!\!_{j_1}$}}{\TC~[tnpos=l]{4}~[tnpos=r]{$\!\!_{j_2}$}}
    \pstree{\TC~[tnpos=l]{3}~[tnpos=r]{$\!\!_{j_1}$}}{\TC~[tnpos=l]{5}~[tnpos=r]{$\!\!_{j_2}$}}}
    }
    \end{tabular}
\caption{Some trees which belong to $LTS(I)$.}
\label{Ito-trees-in-LTS(L0):Wm}
\end{center}
\end{figure}

For the expansion of the expectation of some functional of the
diffusion process $(X_t)_{t \in I}$ which is solution of the
$d$-dimensional It{\^o} SDE
\begin{equation} \label{Ito-SDE1-autonom-Wm}
    X_t = x_0 + \int_{t_0}^t a(X_s) \, ds + \sum_{j=1}^m \int_{t_0}^t
    b^j(X_s) \, dW_s^j
\end{equation}
with constant initial value $x_0 \in \mathbb{R}^d$ we have to
introduce the subset $LTS(I)$ of $LTS$.
\begin{Def} \label{Ito-subset-LTS(I)}
    Let $LTS(I)$ denote the set of trees $\textbf{t}=(\textbf{t}',\textbf{t}'') \in LTS$ having
    as
    root $\gamma = \,\,$~\pstree[treemode=U, dotstyle=otimes, dotsize=3.2mm, levelsep=0.1cm, radius=1.6mm, treefit=loose]
    {\Tn}{
    \pstree[treemode=U, dotstyle=otimes, dotsize=3.2mm, levelsep=0cm, radius=1.6mm, treefit=loose]
    {\Tdot} {}
    }$\,\,\,$
    and which can be constructed by a finite number
    of steps of the form
    \begin{enumerate}[a)]
        \item adding a deterministic node $\tau = $
        \pstree[treemode=U, dotstyle=otimes, dotsize=3.2mm, levelsep=0.1cm, radius=1.6mm, treefit=loose]
        {\Tn}{
        \pstree[treemode=U, dotstyle=otimes, dotsize=3.2mm, levelsep=0cm, radius=1.6mm, treefit=loose]
        {\TC*} {}
        }, or
        \item adding two stochastic nodes $\sigma_{j} = $
        \pstree[treemode=U, dotstyle=otimes, dotsize=3.2mm, levelsep=0.1cm, radius=1.6mm, treefit=loose]
        {\Tn}{
        \pstree[treemode=U, dotstyle=otimes, dotsize=3.2mm, levelsep=0cm, radius=1.6mm, treefit=loose]
        {\TC~[tnpos=r]{$\!\!_{j}$}} {}
        }, where both have the same variable index $j \in J$ and neither of them is father of the other.
    \end{enumerate}
    The nodes have to be labelled in the same order as they have been
    added by the construction of the tree. Further $TS(I) = LTS(I)/ \sim$
    denotes the equivalence class under the relation of
    Definition~\ref{St-tree-equivalence:Wm} restricted to $LTS(I)$
    and $\alpha_I(\textbf{t})$ denotes the cardinality of $\textbf{t}$ in $LTS(I)$.
\end{Def}
Since the number of stochastic nodes is always even, the order
$\rho(\textbf{t})$ has to be an integer and $\textbf{t}$ owns the
variable indices $j_1, \ldots, j_{s(\textbf{t})/2}$. See
Figure~\ref{Ito-trees-in-LTS(L0):Wm} for some examples of $TS(I)$
up to order two.
\noindent Now the first main result providing a compact way of
representing the expectation of some functional of the solution of
an It{\^o} stochastic differential equation can be stated.
\begin{The} \label{Ito-tree-expansion-exact-sol:Wm}
    Let $(X_t)_{t \in I}$ be the solution of the It{\^o}
    stochastic differential equation
    system~(\ref{Ito-SDE1-autonom-Wm})
    with initial value $X_{t_0} = x_0 \in \mathbb{R}^d$.
    Then for $n \in \mathbb{N}_0$ and $f, a^i, b^{i,j} \in C_P^{2(n+1)}(\mathbb{R}^d, \mathbb{R})$
    for $i=1, \ldots, d$ and $j=1, \ldots, m$ we get the
    following:
    \begin{enumerate}
    \item \label{Theorem-Ito-part1:Wm} For every $t \in [t_0,T]$ the following truncated
    expansion holds:
    \begin{equation}
    \begin{split} \label{Ito-tree-expansion-exact-sol-formula1:Wm}
        E^{t_0,x_0}(f(X_{t})) = &\sum_{\substack{\textbf{t} \in LTS(I) \\
        \rho(\textbf{t}) \leq n}} \,\,
        \sum_{j_1, \ldots, j_{s(\textbf{t})/2}=1}^m
        \frac{F(\textbf{t})(x_0) }{2^{s(\textbf{t})/2} \, \rho(\textbf{t})!} \, (t-t_0)^{\rho(\textbf{t})} +
        \mathcal{R}_n(t,t_0) \\ %\label{Ito-tree-expansion-exact-sol-formula2:Wm}
        = &\sum_{\substack{\textbf{t} \in TS(I) \\ \rho(\textbf{t}) \leq n}} \,\,
        \sum_{j_1, \ldots, j_{s(\textbf{t})/2}=1}^m
        \frac{\alpha_I(\textbf{t}) \, F(\textbf{t})(x_0) }{2^{s(\textbf{t})/2}
        \, \rho(\textbf{t})!} \, (t-t_0)^{\rho(\textbf{t})} +
        \mathcal{R}_n(t,t_0)
    \end{split}
    \end{equation}
    \item \label{Theorem-Ito-part2:Wm} There exists a
    random variable $\xi : \Omega \rightarrow \, ]t_0,t[$ so that holds:
    \begin{equation} \label{Ito-tree-expansion-exact-sol-formula2:Wm}
    \begin{split}
        \mathcal{R}_n(t,t_0) = &\sum_{\substack{\textbf{t} \in LTS(I) \\
        \rho(\textbf{t}) = n+1}} \,\,
        \sum_{j_1, \ldots, j_{s(\textbf{t})/2}=1}^m
        \frac{E^{t_0,x_0} (F(\textbf{t})(X_{\xi})) }
        {2^{s(\textbf{t})/2} \, \rho(\textbf{t})!} \, (t-t_0)^{\rho(\textbf{t})} \\
        = &\sum_{\substack{\textbf{t} \in TS(I) \\ \rho(\textbf{t}) = n+1}} \,\,
        \sum_{j_1, \ldots, j_{s(\textbf{t})/2}=1}^m
        \frac{\alpha_I(\textbf{t}) \, E^{t_0,x_0}(F(\textbf{t})(X_{\xi}))
        }{2^{s(\textbf{t})/2} \, \rho(\textbf{t})!} \,
        (t-t_0)^{\rho(\textbf{t})}
    \end{split}
    \end{equation}
    \end{enumerate}
\end{The}
{\bf{Proof.}} Consider the diffusion operator\footnote{$L^0$ is
also called the generator of the solution $X$ (see, e.g.,
\cite{Dy65a,Ok98,SV79}).} $L^0$ and the operator $L^j$ defined
by~(\ref{Operator_L0_L1}). Then, by reapplication of It{\^o}'s
formula as in~(\ref{Taylor-recursive-intro:1}) we get
\begin{equation}
    E^{t_0,x_0}(f(X_t))
    = \sum_{i=0}^n \frac{(L^0)^i f(x_0)}{i!} \, (t-t_0)^i + \mathcal{R}_n^*(t,t_0)
\end{equation}
with
\begin{equation} \label{Ito-theorem-proof-eq001}
    \mathcal{R}_n^*(t,t_0) =
    E^{t_0,x_0} \big( \int_{t_0}^t \int_{t_0}^{s_n} \ldots \int_{t_0}^{s_1}
    (L^0)^{n+1} f(X_{s}) \, ds \, ds_1 \ldots ds_n \big).
\end{equation}
Here, it should be remarked that the first moment of a multiple
It{\^o} integral vanishes if it has at least one integration
w.r.t.\ a Wiener process (see~\cite{KP99}).
Thus, in order to prove
(\ref{Ito-tree-expansion-exact-sol-formula1:Wm}) we have to show
the following equation
\begin{equation}
    \sum_{i=0}^{n} \frac{(L^0)^i f(x_0)}{i!} \, (t-t_0)^i =
    \sum_{\substack{\textbf{t} \in LTS(I) \\ \rho(\textbf{t}) \leq n}} \,\,
    \sum_{j_1, \ldots, j_{s(\textbf{t})/2}=1}^m \frac{F(\textbf{t})(x_0)}
    {2^{s(\textbf{t})/2} \, \rho(\textbf{t})!} \, (t-t_0)^{\rho(\textbf{t})} .
\end{equation}
Therefore, for every $n \in \mathbb{N}_0$ we finally have to prove
\begin{equation} \label{Ito-proof-to-show-1}
    (L^0)^n f(x_0) =
    \sum_{\substack{\textbf{t} \in LTS(I) \\ \rho(\textbf{t}) = n}} \,\,
    \sum_{j_1, \ldots, j_{s(\textbf{t})/2}=1}^m
    \frac{F(\textbf{t})(x_0)}{2^{s(\textbf{t})/2}} .
\end{equation}
We proceed by induction. The case $n=0$ is trivial: There is only
one element $\gamma = \,\,$
    \pstree[treemode=U, dotstyle=otimes, dotsize=3.2mm, levelsep=0.1cm, radius=1.6mm, treefit=loose]
    {\Tn}{
    \pstree[treemode=U, dotstyle=otimes, dotsize=3.2mm, levelsep=0cm, radius=1.6mm, treefit=loose]
    {\Tdot} {}
    }
    $\,$
in $LTS(I)$ with $\rho(\gamma)=0$ and we get with $(L^0)^0 \equiv
\Id$ that $(L^0)^0 f(x_0) = f(x_0) = F(\gamma)(x_0)$. The step
$n=1$ is performed as well for better understanding. In this case
two different trees $\textbf{t}_{1.1} = (\tau)$ and
$\textbf{t}_{1.2} = (\sigma_{j_1}, \sigma_{j_1})$ in $LTS(I)$,
both of order 1, have to be considered. For these two trees we
have
\begin{alignat}{2}
    L^0 \left( F(\gamma)(x_0) \right)
    = & \left( \sum_{k=1}^d a^k \frac{\partial f}{\partial
    x^k} + \frac{1}{2} \sum_{k,l=1}^d
    \sum_{j=1}^m b^{k,j} \, b^{l,j} \, \frac{\partial^2
    f}{\partial x^k \partial x^l} \right) (x_0) \notag \\
    = & \,\, F(\textbf{t}_{1.1})(x_0) + \frac{1}{2}
    \sum_{j_1=1}^m F(\textbf{t}_{1.2})(x_0) \notag \\
    = & \sum_{\substack{\textbf{t} \in LTS(I) \\ \rho(\textbf{t})=1}} \,\,
    \sum_{j_1, \ldots, j_{s(\textbf{t})/2}=1}^m
    \frac{F(\textbf{t})(x_0)}{2^{s(\textbf{t})/2}}. \notag
\end{alignat}
Under the assumption that equation (\ref{Ito-proof-to-show-1})
holds for $n \in \mathbb{N}_0$ we now proceed to give a prove for
the case $n+1$.
Therefore, we get
\begin{equation} \label{Ito-proof-eqn9}
\begin{split}
    (L^0)^{n+1} f(x_0) = & \, (L^0) (L^0)^n f(x_0) =
    (L^0) \sum_{\substack{\textbf{t} \in LTS(I) \\
    \rho(\textbf{t})=n}} \,\, \sum_{j_1, \ldots, j_{s(\textbf{t})/2}=1}^m
    \frac{F(\textbf{t})(x_0)}{2^{s(\textbf{t})/2}} \\
    = & \sum_{\substack{\textbf{t} \in LTS(I) \\ \rho(\textbf{t})=n}} \,\,
    \sum_{j_1, \ldots, j_{s(\textbf{t})/2}=1}^m
    \frac{1}{2^{s(\textbf{t})/2}} \,
    \Big( \sum_{k=1}^d a^k \, \frac{\partial}{\partial
    x^k} F(\textbf{t}) \\
    & + \frac{1}{2} \sum_{k,l=1}^d
    \sum_{j=1}^m b^{k,j} \, b^{l,j} \,
    \frac{\partial^2}{\partial x^k \partial x^l}
    F(\textbf{t}) \Big)(x_0).
\end{split}
\end{equation}
Next consider the first part of the operator $L^0$ and apply
Lemma~\ref{Lem-tree-growth}. Then for $\textbf{u} \in LTS(I)$ with
$\rho(\textbf{u})=n$ and $l(\textbf{u})=\lambda$ we get
\begin{equation} \label{Ito-proof-eqn7}
    \sum_{k=1}^d a^k(x_0) \frac{\partial}{\partial x^k}
    \frac{F(\textbf{u})(x_0)}{2^{s(\textbf{u})/2}}
    = \sum_{\textbf{t} \in H_1(\textbf{u})}
    \frac{F(\textbf{t})(x_0)}{2^{s(\textbf{t})/2}}
\end{equation}
where $H_1(\textbf{u})$ is the set of all trees $\textbf{t} \in
LTS(I)$ with $\rho(\textbf{t})=n+1$, $\textbf{t}' |_{\{2, \ldots,
\lambda\} } = \textbf{u}'$, $\textbf{t}'' |_{\{1, \ldots,
\lambda\} } = \textbf{u}''$ and $\textbf{t}''(\lambda+1) = \tau$.
Clearly $s(\textbf{u})=s(\textbf{t})$ holds for all
$\textbf{t} \in H_1(\textbf{u})$.\\ \\
Now we proceed by considering the second part of the operator
$L^0$. Having Lemma~\ref{Lem-tree-growth-2} in mind, for
$\textbf{u} \in LTS(I)$ with $\rho(\textbf{u})=n$ and
$l(\textbf{u})=\lambda$ we calculate
\begin{equation} \label{Ito-proof-eqn6}
\begin{split}
    & \frac{1}{2} \sum_{k,l=1}^d \sum_{j=1}^m
    b^{k,j}(x_0) \, b^{l,j}(x_0) \,
    \frac{\partial^2}{\partial x^k \partial x^l} \,\,
    \sum_{j_1, \ldots, j_{s(\textbf{u})/2}=1}^m
    \frac{F(\textbf{u})(x_0)}{2^{s(\textbf{u})/2}} \\
    & = \sum_{\textbf{t} \in H_2(\textbf{u})} \,\,
    \sum_{j_1, \ldots, j_{s(\textbf{t})/2}=1}^m
    \frac{F(\textbf{t})(x_0)}{2^{s(\textbf{t})/2}}
\end{split}
\end{equation}
where $H_2(\textbf{u})$ is the set of all trees $\textbf{t} =
(\textbf{t}', \textbf{t}'') \in LTS(I)$ with $l(\textbf{t}) =
\lambda + 2$ nodes, of order $\rho(\textbf{t}) = n+1$, with
$\textbf{t}' |_{\{2,\ldots,\lambda\}} = \textbf{u}'$,
$\textbf{t}'' |_{\{1,\ldots,\lambda\}} = \textbf{u}''$ and
$\textbf{t}''(\lambda+1) = \textbf{t}''(\lambda+2) = \sigma_j =
\sigma_{j_{s(\textbf{u})/2+1}} = \sigma_{j_{s(\textbf{t})/2}}$
since $s(\textbf{t}) = s(\textbf{u}) + 2$
for all $\textbf{t} \in H_2(\textbf{u})$. \\ \\
Taking together the results (\ref{Ito-proof-eqn7}) and
(\ref{Ito-proof-eqn6}), the following equation
\begin{equation} \label{Ito-proof-eqn10}
\begin{split}
    (L^0) \, \sum_{j_1, \ldots, j_{s(\textbf{u})/2} =1}^m
    \frac{F(\textbf{u})(x_0)}{2^{s(\textbf{u})/2}} = &\sum_{\textbf{t} \in H_1(\textbf{u})} \,\,
    \sum_{j_1, \ldots, j_{s(\textbf{t})/2}=1}^m
    \frac{F(\textbf{t})(x_0)}{2^{s(\textbf{t})/2}} \\
    & + \sum_{\textbf{t} \in H_2(\textbf{u})} \,\,
    \sum_{j_1, \ldots, j_{s(\textbf{t})/2}=1}^m
    \frac{F(\textbf{t})(x_0)}{2^{s(\textbf{t})/2}} \\
    = &\sum_{\textbf{t} \in H_1(\textbf{u}) \cup H_2(\textbf{u})} \,\,
    \sum_{j_1, \ldots, j_{s(\textbf{t})/2}=1}^m
    \frac{F(\textbf{t})(x_0)}{2^{s(\textbf{t})/2}}
\end{split}
\end{equation}
holds for every $\textbf{u} \in LTS(I)$ with $\rho(\textbf{u})=n$.
Now it is easily seen that
\begin{equation} \label{Ito-proof-eqn11}
    \bigcup_{\substack{\textbf{u} \in LTS(I) \\ \rho(\textbf{u})=n}}
    \left( H_1(\textbf{u}) \cup H_2(\textbf{u}) \right) = \{ \textbf{t} \in LTS(I) :
    \rho(\textbf{t}) = n+1\} .
\end{equation}
By applying (\ref{Ito-proof-eqn10}) and (\ref{Ito-proof-eqn11}) to
equation (\ref{Ito-proof-eqn9}) we arrive  at
\begin{equation}
\begin{split}
    & (L^0) \sum_{\substack{\textbf{t} \in LTS(I) \\ \rho(\textbf{t})=n}} \,\,
    \sum_{j_1, \ldots, j_{s(\textbf{t})/2}=1}^m
    \frac{F(\textbf{t})(x_0)}{2^{s(\textbf{t})/2}}
    = \sum_{\substack{\textbf{t} \in LTS(I) \\ \rho(\textbf{t})=n+1}} \,\,
    \sum_{j_1, \ldots, j_{s(\textbf{t})/2}=1}^m
    \frac{F(\textbf{t})(x_0)}{2^{s(\textbf{t})/2}}
\end{split}
\end{equation}
which completes the proof of the first part of
Theorem~\ref{Ito-tree-expansion-exact-sol:Wm}. \\ \\
Finally, we have to prove that $\mathcal{R}_n(t,t_0) =
\mathcal{R}_n^*(t,t_0)$. Due to the Existence and Uniqueness
Theorem, there exists $N \subset \Omega$ with $P(N) = 0$ such that
for all $\omega \in \Omega \setminus N$ the solution $X_t(\omega)$
of the Ito SDE~(\ref{Intro-Ito-St-SDE1-integralform-Wm}) is a
continuous function of $t$. Now $f, a^i, b^{i,j} \in
C_P^{2(n+1)}(\mathbb{R}^d,\mathbb{R})$ for $i=1, \ldots, d$ and
$j=1, \ldots, m$ implies that $(L^0)^{n+1} f \in C_P(\mathbb{R}^d,
\mathbb{R})$. Then for every $\omega \in \Omega \setminus N$ there
exists $\xi(\omega) \in \, ]t_0,t[$ so that
\begin{equation*}
    \begin{split}
    \int_{t_0}^t \int_{t_0}^{s_n} \ldots \int_{t_0}^{s_1}
    (L^0)^{n+1} f(X_{s}(\omega)) \, ds \, ds_1 \ldots ds_n =
    \frac{(L^0)^{n+1} f(X_{\xi}(\omega))}{(n+1)!} \,
    (t-t_0)^{n+1}
    \end{split}
\end{equation*}
holds. Applying (\ref{Ito-proof-to-show-1}) and taking the
expectation completes the proof. \hfill $\Box$ \\ \\
A spreadsheet containing all trees $\textbf{t}$ of the set $TS(I)$
up to order 2.0 together with the corresponding graph can be found
in Table~\ref{tabelle1}.
\section{Truncated Taylor expansions for Stratonovich SDE Systems}
\label{rooted-tree-analysis-stratonovich-SDE-systems:Wm}
\begin{figure}[H]
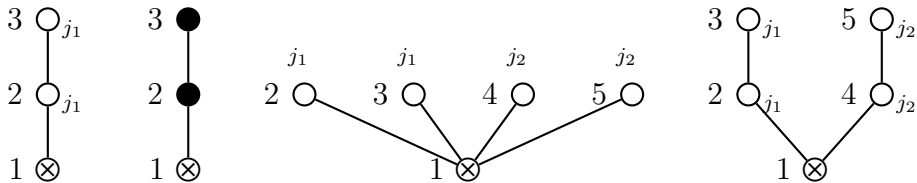

\begin{center}
    \begin{tabular}{ccccccccc}
    \pstree[treemode=U, dotstyle=otimes, dotsize=3.2mm, levelsep=0.1cm, radius=1.6mm, treefit=loose]
    {\Tn}{
    \pstree[treemode=U, dotstyle=otimes, dotsize=3.2mm, levelsep=1cm, radius=1.6mm, treefit=loose]
    {\Tdot~[tnpos=l]{1 }} {\pstree{\TC~[tnpos=l]{2}~[tnpos=r]{$\!\!_{j_1}$}}{\TC~[tnpos=l]{3}~[tnpos=r]{$\!\!_{j_1}$}}}
    }
    & \qquad &
    \pstree[treemode=U, dotstyle=otimes, dotsize=3.2mm, levelsep=0.1cm, radius=1.6mm, treefit=loose]
    {\Tn}{
    \pstree[treemode=U, dotstyle=otimes, dotsize=3.2mm, levelsep=1cm, radius=1.6mm, treefit=loose]
    {\Tdot~[tnpos=l]{1 }} {\pstree{\TC*~[tnpos=l]{2}}{\TC*~[tnpos=l]{3}}}
    }
    & \qquad &
    \pstree[treemode=U, dotstyle=otimes, dotsize=3.2mm, levelsep=0.1cm, radius=1.6mm, treefit=loose]
    {\Tn}{
    \pstree[treemode=U, dotstyle=otimes, dotsize=3.2mm, levelsep=1cm, radius=1.6mm, treefit=loose]
    {\Tdot~[tnpos=l]{1 }} {\TC~[tnpos=l]{2}~[tnpos=a]{$\!\!_{j_1}$} \TC~[tnpos=l]{3}~[tnpos=a]{$\!\!_{j_1}$}
    \TC~[tnpos=l]{4}~[tnpos=a]{$\!\!_{j_2}$} \TC~[tnpos=l]{5}~[tnpos=a]{$\!\!_{j_2}$}}
    }
    & \quad &
    \pstree[treemode=U, dotstyle=otimes, dotsize=3.2mm, levelsep=0.1cm, radius=1.6mm, treefit=loose]
    {\Tn}{
    \pstree[treemode=U, dotstyle=otimes, dotsize=3.2mm, levelsep=1cm, radius=1.6mm, treefit=loose]
    {\Tdot~[tnpos=l]{1 }} {\pstree{\TC~[tnpos=l]{2}~[tnpos=r]{$\!\!_{j_1}$}}{\TC~[tnpos=l]{3}~[tnpos=r]{$\!\!_{j_1}$}}
    \pstree{\TC~[tnpos=l]{4}~[tnpos=r]{$\!\!_{j_2}$}}{\TC~[tnpos=l]{5}~[tnpos=r]{$\!\!_{j_2}$}}}
    }
    \end{tabular}
\caption{Some trees which belong to $TS(S)$.}
\label{St-trees-in-LTS(L0):Wm}
\end{center}
\end{figure}

An expansion of the expectation of some functional applied to the
solution of a Stratonovich SDE system can be stated in a similar
way as introduced for It{\^o} SDE systems. Therefore we consider
the following SDE w.r.t.\ Stratonovich calculus
\begin{equation} \label{St-SDE1-autonom-Wm}
    X_t = x_0 + \int_{t_0}^t \underline{a}(X_s) \, ds
    + \sum_{j=1}^m \int_{t_0}^t b^j(X_s) \, \circ dW_s .
\end{equation}
The solution $(X_t)_{t \in I}$ of the Stratonovich SDE is also
solution of a corresponding It{\^o} SDE as
in~(\ref{Ito-SDE1-autonom-Wm}) and therefore also a diffusion
process, however with modified drift
\begin{equation} \label{St-SDE1-Ito-correction-tilde:a}
    \tilde{a}^i(x) = \underline{a}^i(x) + \frac{1}{2} \sum_{k=1}^d \sum_{l=1}^m b^{k,l}(x)
    \, \frac{\partial b^{i,l}}{\partial x^k}(x)
\end{equation}
for $i=1, \ldots, d$. Again we introduce a special subset $LTS(S)$
of $LTS$.
\begin{Def} \label{St-subset-LTS(S)}
    Let $LTS(S)$ denote the set of trees $\textbf{t}=(\textbf{t}',\textbf{t}'') \in LTS$ having
    as
    root $\gamma = \,\,$~\pstree[treemode=U, dotstyle=otimes, dotsize=3.2mm, levelsep=0.1cm, radius=1.6mm, treefit=loose]
    {\Tn}{
    \pstree[treemode=U, dotstyle=otimes, dotsize=3.2mm, levelsep=0cm, radius=1.6mm, treefit=loose]
    {\Tdot} {}
    }$\,\,\,$
    and which can be constructed by a finite number
    of steps of the form
    \begin{enumerate}[a)]
        \item adding a deterministic node $\tau = $
        \pstree[treemode=U, dotstyle=otimes, dotsize=3.2mm, levelsep=0.1cm, radius=1.6mm, treefit=loose]
        {\Tn}{
        \pstree[treemode=U, dotstyle=otimes, dotsize=3.2mm, levelsep=0cm, radius=1.6mm, treefit=loose]
        {\TC*} {}
        }, or
        \item adding two stochastic nodes $\sigma_j = $
        \pstree[treemode=U, dotstyle=otimes, dotsize=3.2mm, levelsep=0.1cm, radius=1.6mm, treefit=loose]
        {\Tn}{
        \pstree[treemode=U, dotstyle=otimes, dotsize=3.2mm, levelsep=0cm, radius=1.6mm, treefit=loose]
        {\TC~[tnpos=r]{$\!\!_{j}$}} {}
        }, where both nodes have the same variable index $j \in J$.
    \end{enumerate}
    The nodes have to be labelled in the same order as they have been
    added by the construction of the tree. Further $TS(S) = LTS(S)/ \sim$
    denotes the equivalence class under the relation of
    Definition~\ref{St-tree-equivalence:Wm} restricted to $LTS(S)$
    and $\alpha_S(\textbf{t})$ denotes the cardinality of $\textbf{t}$ in $LTS(S)$.
\end{Def}
The construction of the trees $\textbf{t} \in LTS(S)$ forces the
number of stochastic nodes to be even, the order
$\rho(\textbf{t})$ has to be an integer and $\textbf{t}$ owns the
variable indices $j_1, \ldots, j_{s(\textbf{t})/2}$.
Figure~\ref{St-trees-in-LTS(L0):Wm} presents some examples up to
order two of the set $TS(S)$. As the construction of the trees in
$LTS(I)$ is more restrictive than the ones in $LTS(S)$, we have
$LTS(I) \subset LTS(S)$. As an example, the first and the last
tree of
Figure~\ref{St-trees-in-LTS(L0):Wm} don't belong to $TS(I)$. \\
\noindent The following result can be proved in nearly the same
way as Theorem~\ref{Ito-tree-expansion-exact-sol:Wm}.
\begin{The} \label{St-tree-expansion-exact-sol:Wm}
    Let $(X_t)_{t \in I}$ be the solution of the Stratonovich
    stochastic differential equation
    system~(\ref{St-SDE1-autonom-Wm})
    with initial value $X_{t_0} = x_0 \in \mathbb{R}^d$.
    Then for $n \in \mathbb{N}_0$ and $f, \tilde{a}^i, b^{i,j} \in C_P^{2(n+1)}(\mathbb{R}^d, \mathbb{R})$
    for $i=1, \ldots, d$ and $j=1, \ldots, m$ we get the
    following:
    \begin{enumerate}
    \item \label{Theorem-St-part1:Wm} For every $t \in [t_0,T]$ the following truncated
    expansion holds:
    \begin{equation}
    \begin{split} \label{Strato-tree-expansion-exact-sol-formula1:Wm}
        E^{t_0,x_0}(f(X_{t})) = &\sum_{\substack{\textbf{t} \in LTS(S) \\
        \rho(\textbf{t}) \leq n}} \,\,
        \sum_{j_1, \ldots, j_{s(\textbf{t})/2}=1}^m
        \frac{F(\textbf{t})(x_0) }{2^{s(\textbf{t})/2} \, \rho(\textbf{t})!} \, (t-t_0)^{\rho(\textbf{t})} +
        \underline{\mathcal{R}}_n(t,t_0) \\
        = &\sum_{\substack{\textbf{t} \in TS(S) \\ \rho(\textbf{t}) \leq n}} \,\,
        \sum_{j_1, \ldots, j_{s(\textbf{t})/2}=1}^m
        \frac{\alpha_S(\textbf{t}) \, F(\textbf{t})(x_0) }{2^{s(\textbf{t})/2}
        \, \rho(\textbf{t})!} \, (t-t_0)^{\rho(\textbf{t})} +
        \underline{\mathcal{R}}_n(t,t_0)
    \end{split}
    \end{equation}
    \item \label{Theorem-St-part2:Wm} There exists a
    random variable $\xi : \Omega \rightarrow \, ]t_0,t[$ so that holds:
    \begin{equation} \label{Strato-tree-expansion-exact-sol-formula2:Wm}
    \begin{split}
        \underline{\mathcal{R}}_n(t,t_0) = &\sum_{\substack{\textbf{t} \in LTS(S) \\
        \rho(\textbf{t}) = n+1}} \,\,
        \sum_{j_1, \ldots, j_{s(\textbf{t})/2}=1}^m
        \frac{E^{t_0,x_0} (F(\textbf{t})(X_{\xi})) }
        {2^{s(\textbf{t})/2} \, \rho(\textbf{t})!} \, (t-t_0)^{\rho(\textbf{t})} \\
        = &\sum_{\substack{\textbf{t} \in TS(S) \\ \rho(\textbf{t}) = n+1}} \,\,
        \sum_{j_1, \ldots, j_{s(\textbf{t})/2}=1}^m
        \frac{\alpha_S(\textbf{t}) \, E^{t_0,x_0}(F(\textbf{t})(X_{\xi}))
        }{2^{s(\textbf{t})/2} \, \rho(\textbf{t})!} \,
        (t-t_0)^{\rho(\textbf{t})}
    \end{split}
    \end{equation}
    \end{enumerate}
\end{The}
{\bf{Proof.}} The proof is very similar to the one of
Theorem~\ref{Ito-tree-expansion-exact-sol:Wm}. Since the solution
$(X_t)_{t \in I}$ of the autonomous Stratonovich stochastic
differential equation~(\ref{St-SDE1-autonom-Wm}) is also solution
of a corresponding It{\^o} stochastic differential equation with
the modified drift
$\tilde{a}$~(\ref{St-SDE1-Ito-correction-tilde:a}), we refer to
the corresponding It{\^o} stochastic differential equation in the
following.
Then the diffusion operator $L^0$ results in
\begin{equation}
    \begin{split}
    L^0 = &\sum_{k=1}^d \tilde{a}^k \, \frac{\partial}{\partial x^k} +
    \frac{1}{2} \sum_{k,l=1}^d \sum_{j=1}^m b^{k,j}  b^{l,j} \,
    \frac{\partial^2}{\partial x^k \partial x^l} \\
    = & \sum_{k=1}^d \underline{a}^k \, \frac{\partial}{\partial x^k}
    + \sum_{j=1}^m \left( \frac{1}{\sqrt{2}} \sum_{k=1}^d b^{k,j}
    \, \frac{\partial}{\partial x^k} \right)^2
    \end{split}
\end{equation}
while the operator $L^j$ remains as defined
in~(\ref{Operator_L0_L1}). By reapplication of It{\^o}'s formula
as performed in~(\ref{Taylor-recursive-intro:1}) we get
\begin{equation}
    E^{t_0,x_0}(f(X_t)) = \sum_{i=0}^n \frac{(L^0)^i f(x_0)}{i!}
    (t-t_0)^i + \underline{\mathcal{R}}^*_n(t,t_0)
\end{equation}
with $\underline{\mathcal{R}}_n^*(t,t_0)$ in a similar way as in
(\ref{Ito-theorem-proof-eq001}).
Thus, to prove (\ref{Strato-tree-expansion-exact-sol-formula1:Wm})
it is sufficient to show for every $n \in \mathbb{N}_0$ the
equation
\begin{equation} \label{St-proof-to-show-1}
    (L^0)^n f(x_0) =
    \sum_{\substack{\textbf{t} \in LTS(S) \\ \rho(\textbf{t}) = n}} \,\,
    \sum_{j_1, \ldots, j_{s(\textbf{t})/2}=1}^m
    \frac{F(\textbf{t})(x_0)}{2^{s(\textbf{t})/2}} .
\end{equation}
The proof proceeds by induction. The case $n=0$ is trivial. Again
step $n=1$ is performed for better understanding. In this case
three different trees $\textbf{t}_{1.1}=(\tau)$,
$\textbf{t}_{1.2}=(\sigma_{j_1}, \sigma_{j_1})$ and
$\textbf{t}_{1.3}=(\{ \sigma_{j_1} \}_{j_1})$ in $LTS(S)$, all of
order 1, have to be considered.
For these three trees we have
\begin{alignat}{2}
    L^0 \left( F(\gamma)(x_0) \right) %\notag \\
    & = \,\, F(\textbf{t}_{1.1})(x_0) + \tfrac{1}{2} \sum_{j_1=1}^m F(\textbf{t}_{1.3})(x_0) +
    \tfrac{1}{2} \sum_{j_1=1}^m F(\textbf{t}_{1.2})(x_0) \notag \\
    & = \sum_{\substack{\textbf{t} \in LTS(S) \\ \rho(\textbf{t}) = 1}} \,\,
    \sum_{j_1, \ldots, j_{s(\textbf{t})/2}=1}^m
    \frac{F(\textbf{t})(x_0)}{2^{s(\textbf{t})/2}}. \notag
\end{alignat}
Under the assumption that equation (\ref{St-proof-to-show-1})
holds for $n \in \mathbb{N}_0$ we proceed to prove the case $n+1$.
Therefore, we get
\begin{equation} \label{St-proof-eqn9}
\begin{split}
    (L^0)^{n+1} f(x_0)
    = & \sum_{\substack{\textbf{t} \in LTS(S) \\ \rho(\textbf{t})=n}} \,\,
    \sum_{j_1, \ldots, j_{s(\textbf{t})/2}=1}^m
    \frac{1}{2^{s(\textbf{t})/2}} \,
    \Big( \sum_{k=1}^d \underline{a}^k \frac{\partial}{\partial
    x^k} F(\textbf{t}) \\
    & + \sum_{j=1}^m \Big( \frac{1}{\sqrt{2}} \sum_{k=1}^d
    b^{k,j} \frac{\partial}{\partial x^k} \Big)^2
    F(\textbf{t}) \Big) (x_0).
\end{split}
\end{equation}
Firstly, we apply Lemma~\ref{Lem-tree-growth} to the first part of
$L^0$ in the same way as in the proof of
Theorem~\ref{Ito-tree-expansion-exact-sol:Wm}, however we denote
for $\textbf{u} \in LTS(S)$ with $\rho(\textbf{u})=n$ and
$l(\textbf{u})=\lambda$ the set $H_1(\textbf{u})$ by
$D(\textbf{u})$ in the following.
Then we proceed by considering the second part of the operator
$L^0$, which is applied in two steps. Keeping
Lemma~\ref{Lem-tree-growth} in mind, as a first step for
$\textbf{u} \in LTS(S)$ with $\rho(\textbf{u})=n$ and
$l(\textbf{u})=\lambda$, we calculate
\begin{equation} \label{St-proof-eqn5}
    \frac{1}{\sqrt{2}} \sum_{j=1}^m
    \sum_{k=1}^d b^{k,j}(x_0) \frac{\partial}{\partial x^k}
    \sum_{j_1, \ldots, j_{s(\textbf{u})/2}=1}^m
    \frac{F(\textbf{u})(x_0)}{2^{s(\textbf{u})/2}}
    = \sum_{\textbf{t} \in S(\textbf{u})} \,\,
    \sum_{j_1, \ldots, j_{s(\textbf{u})/2+1}=1}^m
    \frac{F(\textbf{t})(x_0)}{2^{s(\textbf{t})/2}}
\end{equation}
where $S(\textbf{u})$ denotes the set of trees $\textbf{t} \in
LTS$ with $\rho(\textbf{t})=n+\tfrac{1}{2}$, $\textbf{t}'
|_{\{2,\ldots,\lambda\}} = \textbf{u}'$, $\textbf{t}'' |_{\{1,
\ldots, \lambda\}} = \textbf{u}''$ and $\textbf{t}''(\lambda + 1)
= \sigma_j = \sigma_{s(\textbf{u})/2+1}$. Here we have
$s(\textbf{t}) = s(\textbf{u}) + 1$ for all $\textbf{t} \in
S(\textbf{u})$. For the second step, we repeat the first step for
$\textbf{u}$ replaced by $\textbf{t}=(\textbf{t}',\textbf{t}'')
\in S(\textbf{u})$. \\ \\
As a result of this we get for $\textbf{u} \in LTS(S)$ with
$\rho(\textbf{u})=n$ and $l(\textbf{u})=\lambda$
\begin{equation} \label{St-proof-eqn6}
\begin{split}
    \sum_{j=1}^m & \left( \frac{1}{\sqrt{2}} \sum_{k=1}^d b^{k,j}(x_0)
    \frac{\partial}{\partial x^k} \right)^2
    \sum_{j_1, \ldots, j_{s(\textbf{u})/2}=1}^m
    \frac{F(\textbf{u})(x_0)}{2^{s(\textbf{u})/2}} \\
    & =
    \sum_{\textbf{t}_1 \in S(\textbf{u})} \,\,
    \sum_{j_1, \ldots, j_{s(\textbf{u})/2+1}=1}^m
    \frac{1}{\sqrt{2}} \sum_{k=1}^d b^{k,j_{s(\textbf{u})/2+1}}(x_0)
    \frac{\partial}{\partial x^k}
    \frac{F(\textbf{t}_1)(x_0)}{2^{s(\textbf{t}_1)/2}} \\
    & = \sum_{\textbf{t}_1 \in S(\textbf{u})} \,\, \sum_{\textbf{t}_2 \in S(\textbf{t}_1)} \,\,
    \sum_{j_1, \ldots, j_{s(\textbf{t}_2)/2}=1}^m
    \frac{F(\textbf{t}_2)(x_0)}{2^{s(\textbf{t}_2)/2}} \\
    & = \sum_{\textbf{t} \in S(S(\textbf{u}))} \,\,
    \sum_{j_1, \ldots, j_{s(\textbf{t})/2}=1}^m
    \frac{F(\textbf{t})(x_0)}{2^{s(\textbf{t})/2}}
\end{split}
\end{equation}
where $S(S(\textbf{u})) = \bigcup_{\textbf{t}_1 \in S(\textbf{u})}
S(\textbf{t}_1)$ is the set of all trees $\textbf{t} =
(\textbf{t}', \textbf{t}'') \in LTS(S)$ with $l(\textbf{t}) =
\lambda + 2$ nodes, order $\rho(\textbf{t}) = n+1$, with
$\textbf{t}' |_{\{2,\ldots,\lambda\}} = \textbf{u}'$,
$\textbf{t}'' |_{\{1,\ldots,\lambda\}} = \textbf{u}''$ and
$\textbf{t}''(\lambda+1) = \textbf{t}''(\lambda+2) = \sigma_j =
\sigma_{j_{s(\textbf{u})/2+1}} = \sigma_{j_{s(\textbf{t})/2}}$.
Again we have $s(\textbf{t}) = s(\textbf{u}) + 2$
for all $\textbf{t} \in S(S(\textbf{u}))$. \\ \\
Taking together now the results for the first and the second part
of $L^0$, the equation
\begin{equation} \label{St-proof-eqn10}
\begin{split}
    (L^0) \sum_{j_1, \ldots, j_{s(\textbf{u})/2}=1}^m
    \frac{F(\textbf{u})(x_0)}{2^{s(\textbf{u})/2}}
    = &\sum_{\textbf{t} \in D(\textbf{u}) \cup S(S(\textbf{u}))} \,\,
    \sum_{j_1, \ldots, j_{s(\textbf{t})/2}=1}^m
    \frac{F(\textbf{t})(x_0)}{2^{s(\textbf{t})/2}}
\end{split}
\end{equation}
holds for every $\textbf{u} \in LTS(S)$ with $\rho(\textbf{u})=n$.
Now it is easily seen that
\begin{equation} \label{St-proof-eqn11}
    \bigcup_{\substack{\textbf{u} \in LTS(S) \\ \rho(\textbf{u})=n}}
    D(\textbf{u}) \cup S(S(\textbf{u})) = \{ \textbf{t} \in
    LTS(S) : \rho(\textbf{t}) = n+1\}
\end{equation}
and by applying (\ref{St-proof-eqn10}) and (\ref{St-proof-eqn11})
to (\ref{St-proof-eqn9}) we thus arrive at
(\ref{St-proof-to-show-1}) with $n$ replaced by $n+1$, which
completes the proof of the first part of
Theorem~\ref{St-tree-expansion-exact-sol:Wm}.
Finally, we have to prove that $\underline{\mathcal{R}}_n(t,t_0) =
\underline{\mathcal{R}}_n^*(t,t_0)$. This follows by similar
considerations as performed in the proof of
Theorem~\ref{Ito-tree-expansion-exact-sol:Wm}. \hfill $\Box$ \\ \\
A spreadsheet containing all trees $\textbf{t}$ of the set $TS(S)$
up to order 2.0 together with the corresponding graph can be found
in Table~\ref{tabelle1}.
\section{Stochastic Taylor series for It{\^o} and Stratonovich SDE
systems}
\begin{Sat} \label{Proposition-Ito-estimation1:Wm}
For $R_n(t,t_0)=\mathcal{R}_n(t,t_0)$ and for
$R_n(t,t_0)=\underline{\mathcal{R}}_n(t,t_0)$, some $r \in
\mathbb{N}$ and constants $C_1, C_2 > 0$ holds
\begin{equation} \label{Proposition-Ito-estimation1:Wm:formel1}
    | R_n(t,t_0) | \leq C_1 \, (1+ \| x_0
    \|^{2r}) \cdot \exp(C_2 (t-t_0)) \,
    \frac{(t-t_0)^{n+1}}{(n+1)!} .
\end{equation}
\end{Sat}
{\bf{Proof.}} Since $f, a^i, \tilde{a}^i, b^{i,j} \in
C_P^{2(n+1)}(\mathbb{R}^d,\mathbb{R})$ for $i=1, \ldots, d$ and
$j=1, \ldots, m$ we have $(L^0)^{n+1} f \in C_P(\mathbb{R}^d,
\mathbb{R})$. Due to the Existence and Uniqueness Theorem, for
every $T>0$ there exists a constant $C$ which depends only on $T$,
$l \in \mathbb{N}$ and the Lipschitz constant of $a$, $\tilde{a}$
and $b$, so that for all $t \in [t_0,T]$ the estimation
    \begin{equation} \label{esti1}
        E(\| X_t \|^{2l}) \leq (1 + E(\| x_0 \|^{2l})) \cdot \exp(C \,
        (t-t_0))
    \end{equation}
holds \cite{Arn73,KS99}. Therefore, with (\ref{poly_growth}) for
some $r \in \mathbb{N}$ we obtain
\begin{equation}
    \begin{split}
    | R_n(t,t_0) |
    \leq \,\, & \int_{t_0}^t \int_{t_0}^{s_n} \ldots \int_{t_0}^{s_1}
    E^{t_0,x_0} ( | (L^0)^{n+1} f(X_{s}) | ) \, ds \, ds_1 \ldots ds_n \\
    \leq \,\, & \int_{t_0}^t \int_{t_0}^{s_n} \ldots \int_{t_0}^{s_1}
    E^{t_0,x_0} ( K_1 \, (1+ \|X_s \|^{2r}) ) \, ds \, ds_1
    \ldots ds_n \\
    \leq \,\, &\tfrac{(t-t_0)^{n+1}}{(n+1)!} \, K_2 \, (1+
    E^{t_0,x_0}( \| x_0 \|^{2r})) \cdot \exp(K_3 (t-t_0))
    \end{split}
\end{equation}
with some constants $K_1, K_2, K_3 > 0$. \hfill $\Box$
\begin{Kor}
    Let $(X_t)_{t \in I}$ be the solution of the
    stochastic differential equation~(\ref{Intro-Ito-St-SDE1-integralform-Wm})
    with initial value $X_{t_0} = x_0 \in \mathbb{R}^d$.
    Then for $f, a^i, \tilde{a}^i, b^{i,j} \in C_P^{\infty}(\mathbb{R}^d, \mathbb{R})$
    for $i=1, \ldots, d$ and $j=1, \ldots, m$ and for every $t \in [t_0,T]$
    the following expansion holds:
    \begin{equation} \label{Corollar-Ito-eqn1:Wm}
    \begin{split}
        E^{t_0,x_0}(f(X_{t})) = &\sum_{\textbf{t} \in LTS(*)} \,\,
        \sum_{j_1, \ldots, j_{s(\textbf{t})/2}=1}^m
        \frac{F(\textbf{t})(x_0) }{2^{s(\textbf{t})/2} \, \rho(\textbf{t})!} \, (t-t_0)^{\rho(\textbf{t})}
        \\
        = &\sum_{\textbf{t} \in TS(*)} \,\,
        \sum_{j_1, \ldots, j_{s(\textbf{t})/2}=1}^m
        \frac{\alpha_*(\textbf{t}) \, F(\textbf{t})(x_0) }{2^{s(\textbf{t})/2}
        \, \rho(\textbf{t})!} \, (t-t_0)^{\rho(\textbf{t})}
    \end{split}
    \end{equation}
    Here, $*=I$ for the It{\^o} version of SDE~(\ref{Intro-Ito-St-SDE1-integralform-Wm}),
    and $*=S$ for the Stratonovich version of SDE~(\ref{Intro-Ito-St-SDE1-integralform-Wm}).
\end{Kor}
{\bf{Proof.}} This follows from
Theorem~\ref{Ito-tree-expansion-exact-sol:Wm} in case of It{\^o}
calculus and Theorem~\ref{St-tree-expansion-exact-sol:Wm} for
Stratonovich calculus, together with
Proposition~\ref{Proposition-Ito-estimation1:Wm}. For the
convergence of the infinite series~(\ref{Corollar-Ito-eqn1:Wm})
one has to prove that $\lim_{n \rightarrow \infty}
\mathcal{R}_n(t,t_0) = 0$ and $\lim_{n \rightarrow \infty}
\underline{\mathcal{R}}_n(t,t_0) = 0$, which follows
from~(\ref{Proposition-Ito-estimation1:Wm:formel1}). \hfill $\Box$
\section{Example}
Table~\ref{tabelle1} contains all S-trees of the subsets $TS(I)$
and $TS(S)$ up to order two and the corresponding cardinalities
$\alpha_I$ and $\alpha_S$. The cardinalities can be determined
very easily as the number of possibilities to build up the
considered tree due to the corresponding rules of growth in
Definition~\ref{Ito-subset-LTS(I)} and
Definition~\ref{St-subset-LTS(S)}.
\begin{table}[htbp]
\begin{center}
\begin{tabular}[c]{|c|c|c|c|c||c|c|c|c|c|}
    \hline
    $\textbf{t}$ & tree & $\alpha_I$ & $\alpha_S$ & $\rho$ &
    $\textbf{t}$ & tree & $\alpha_I$ & $\alpha_S$ & $\rho$ \\
    \hline
    \hline
    $\textbf{t}_{0.1}$ & $\gamma$ & 1 & 1 & 0 & $\textbf{t}_{2.11}$ & $(\sigma_{j_1},\sigma_{j_1},\sigma_{j_2},\sigma_{j_2})$ & 1 & 1 & 2 \\
    \cline{1-5}
    $\textbf{t}_{1.1}$ & $(\tau)$ & 1 & 1 & 1 & $\textbf{t}_{2.12a}$ & $(\sigma_{j_1},\sigma_{j_1},\{\sigma_{j_2}\}_{j_2})$ & 0 & 2 & 2\\
    $\textbf{t}_{1.2}$ & $(\sigma_{j_1},\sigma_{j_1})$ & 1 & 1 & 1 & $\textbf{t}_{2.12b}$ & $(\sigma_{j_1},\sigma_{j_2},\{\sigma_{j_2}\}_{j_1})$ & 4 & 4 & 2 \\
    $\textbf{t}_{1.3}$ & $(\{\sigma_{j_1}\}_{j_1})$ & 0 & 1 & 1 & $\textbf{t}_{2.13a}$ & $(\sigma_{j_1},\{\sigma_{j_2},\sigma_{j_2}\}_{j_1})$ & 2 & 2 & 2 \\
    \cline{1-5}
    $\textbf{t}_{2.1}$ & $([\tau])$ & 1 & 1 & 2 & $\textbf{t}_{2.13b}$ & $(\sigma_{j_2},\{\sigma_{j_2},\sigma_{j_1}\}_{j_1})$ & 0 & 2 & 2 \\
    $\textbf{t}_{2.2}$ & $(\tau,\tau)$ & 1 & 1 & 2 & $\textbf{t}_{2.14a}$ & $(\sigma_{j_1},\{\{\sigma_{j_2}\}_{j_2}\}_{j_1})$ & 0 & 2 & 2 \\
    $\textbf{t}_{2.3}$ & $([\{\sigma_{j_1}\}_{j_1}])$ & 0 & 1 & 2 & $\textbf{t}_{2.14b}$ & $(\sigma_{j_2},\{\{\sigma_{j_2}\}_{j_1}\}_{j_1})$ & 0 & 2 & 2 \\
    $\textbf{t}_{2.4}$ & $([\sigma_{j_1},\sigma_{j_1}])$ & 1 & 1 & 2 & $\textbf{t}_{2.15a}$ & $(\{\sigma_{j_1}\}_{j_1},\{\sigma_{j_2}\}_{j_2})$ & 0 & 1 & 2 \\
    $\textbf{t}_{2.5}$ & $(\sigma_{j_1},[\sigma_{j_1}])$ & 2 & 2 & 2 & $\textbf{t}_{2.15b}$ & $(\{\sigma_{j_2}\}_{j_1},\{\sigma_{j_2}\}_{j_1})$ & 2 & 2 & 2 \\
    $\textbf{t}_{2.6}$ & $(\{\sigma_{j_1}\}_{j_1},\tau)$ & 0 & 2 & 2 & $\textbf{t}_{2.16}$ & $(\{\sigma_{j_1},\sigma_{j_2},\sigma_{j_2}\}_{j_1})$ & 0 & 1 & 2 \\
    $\textbf{t}_{2.7}$ & $(\sigma_{j_1},\sigma_{j_1},\tau)$ & 2 & 2 & 2 & $\textbf{t}_{2.17a}$ & $(\{\sigma_{j_1},\{\sigma_{j_2}\}_{j_2}\}_{j_1})$ & 0 & 1 & 2 \\
    $\textbf{t}_{2.8}$ & $(\sigma_{j_1},\{\tau\}_{j_1})$ & 2 & 2 & 2 & $\textbf{t}_{2.17b}$ & $(\{\sigma_{j_2},\{\sigma_{j_2}\}_{j_1}\}_{j_1})$ & 0 & 2 & 2 \\
    $\textbf{t}_{2.9}$ & $(\{\{\tau\}_{j_1}\}_{j_1})$ & 0 & 1 & 2 & $\textbf{t}_{2.18}$ & $(\{\{\sigma_{j_2},\sigma_{j_2}\}_{j_1}\}_{j_1})$ & 0 & 1 & 2 \\
    $\textbf{t}_{2.10}$ & $(\{\sigma_{j_1},\tau\}_{j_1})$ & 0 & 1 & 2 & $\textbf{t}_{2.19}$ & $(\{\{\{\sigma_{j_2}\}_{j_2}\}_{j_1}\}_{j_1})$ & 0 & 1 & 2 \\
    \hline
\end{tabular}
\caption{Trees $t \in TS$ for the Taylor expansion of the solution
of a SDE w.r.t.\ either It{\^o} or Stratonovich calculus of order
$0$, $1$ and $2$.} \label{tabelle1}
\end{center}
\end{table}
As an example, we consider the truncated expansion of the
expectation of some functional $f$ of the solution $X_t$ of the
It{\^o} SDE~(\ref{Ito-SDE1-autonom-Wm}). Then, the expansion up to
order two is given by
\begin{equation} \label{Example-allg-expansion}
    \begin{split}
    E^{t_0,x_0} \left(f \left(X_t \right) \right) = &
    F(\textbf{t}_{0.1})(x_0) + F(\textbf{t}_{1.1})(x_0) \cdot
    (t-t_0) + \sum_{j_1=1}^m \frac{F(\textbf{t}_{1.2})(x_0)}{2}
    (t-t_0) \\
    & + \frac{F(\textbf{t}_{2.1})(x_0)}{2} (t-t_0)^2
    + \frac{F(\textbf{t}_{2.2})(x_0)}{2} (t-t_0)^2 \\
    & + \sum_{j_1=1}^m \frac{F(\textbf{t}_{2.4})(x_0)}{4}
    (t-t_0)^2
    + \sum_{j_1=1}^m \frac{F(\textbf{t}_{2.5})(x_0)}{2}
    (t-t_0)^2 \\
    & + \sum_{j_1=1}^m \frac{F(\textbf{t}_{2.7})(x_0)}{2}
    (t-t_0)^2
    + \sum_{j_1=1}^m \frac{F(\textbf{t}_{2.8})(x_0)}{2}
    (t-t_0)^2 \\
    & + \sum_{j_1, j_2 =1}^m \frac{F(\textbf{t}_{2.11})(x_0)}{8}
    (t-t_0)^2
    + \sum_{j_1, j_2 =1}^m \frac{F(\textbf{t}_{2.12b})(x_0)}{2}
    (t-t_0)^2 \\
    & + \sum_{j_1, j_2 =1}^m \frac{F(\textbf{t}_{2.13a})(x_0)}{4}
    (t-t_0)^2
    + \sum_{j_1, j_2 =1}^m \frac{F(\textbf{t}_{2.15b})(x_0)}{4}
    (t-t_0)^2 \\
    & + \mathcal{R}_2(t,t_0) .
    \end{split}
\end{equation}
Finally, we calculate the elementary differentials for the trees
appearing in the truncated
expansion~(\ref{Example-allg-expansion}). As a result of this, we
get
{\allowdisplaybreaks
\begin{eqnarray} % \label{Example-allg-expan-long}
%    \begin{split}
    E^{t_0,x_0} &&(f(X_t)) =
    f(x_0) + \sum_{J=1}^d \frac{\partial f(x_0)}{\partial x^J}
    \, a^J(x_0) \cdot (t-t_0) \notag \\
    && + \frac{1}{2} \sum_{j_1=1}^m \,\, \sum_{J,K=1}^d \frac{\partial^2 f(x_0)}{\partial x^J \partial x^K}
    \, b^{J,j_1}(x_0) \, b^{K,j_1}(x_0) \cdot (t-t_0) \notag \\
    && + \frac{1}{2} \sum_{J,K=1}^d \frac{\partial f(x_0)}{\partial x^J} \,
    \frac{\partial a^J(x_0)}{\partial x^K} \, a^K(x_0) \cdot
    (t-t_0)^2 \notag \\
    && + \frac{1}{2} \sum_{J,K=1}^d
    \frac{\partial^2 f(x_0)}{\partial x^J \partial x^K} \, a^J(x_0) \, a^K(x_0)
    \cdot (t-t_0)^2 \notag \\
    && + \frac{1}{4} \sum_{j_1=1}^m \,\,
    \sum_{J,K,L=1}^d \frac{\partial f(x_0)}{\partial x^J} \,
    \frac{\partial^2 a^J(x_0)}{\partial x^K \partial x^L} \,
    b^{K,j_1}(x_0) \, b^{L,j_1}(x_0) \cdot (t-t_0)^2 \notag \\
    && + \frac{1}{2} \sum_{j_1=1}^m \,\, \sum_{J,K,L=1}^d
    \frac{\partial^2 f(x_0)}{\partial x^J \partial x^K}
    \, b^{J,j_1}(x_0) \, \frac{\partial a^K(x_0)}{\partial x^L}
    \, b^{L,j_1}(x_0) \cdot (t-t_0)^2 \notag \\
    && + \frac{1}{2} \sum_{j_1=1}^m \,\, \sum_{J,K,L=1}^d
    \frac{\partial^3 f(x_0)}{\partial x^J
    \partial x^K \partial x^L} \, b^{J,j_1}(x_0) \, b^{K,j_1}(x_0) \,
    a^L(x_0) \cdot (t-t_0)^2 \notag \\
    && + \frac{1}{2} \sum_{j_1=1}^m \,\, \sum_{J,K,L=1}^d
    \frac{\partial^2 f(x_0)}{\partial x^J \partial x^K} \,
    b^{J,j_1}(x_0) \, \frac{\partial b^{K,j_1}(x_0)}{\partial x^L}
    \, a^L(x_0) \cdot (t-t_0)^2 \notag \\
    && + \frac{1}{8} \sum_{j_1, j_2 =1}^m \,\, \sum_{J,K,L,M=1}^d
    \frac{\partial^4 f(x_0)}{\partial x^J \partial x^K \partial x^L \partial x^M}
    \, b^{J,j_1}(x_0) \, b^{K,j_1}(x_0) \, \times \notag \\
    && \quad \times b^{L,j_2}(x_0) \, b^{M,j_2}(x_0)
    \cdot (t-t_0)^2 \notag \\
    && + \frac{1}{2} \sum_{j_1, j_2 =1}^m \,\, \sum_{J,K,L,M=1}^d
    \frac{\partial^3 f(x_0)}{\partial x^J \partial x^K \partial x^L}
    \, b^{J,j_1}(x_0) \, b^{K,j_2}(x_0) \, \times \notag \\
    && \quad \times \frac{\partial b^{L,j_1}(x_0)}{\partial x^M} \, b^{M,j_2}(x_0)
    \cdot (t-t_0)^2 \notag \\
    && + \frac{1}{4} \sum_{j_1, j_2 =1}^m \,\, \sum_{J,K,L,M=1}^d
    \frac{\partial^2 f(x_0)}{\partial x^J \partial x^K} \, b^{J,j_1}(x_0) \,
    \frac{\partial^2 b^{K,j_1}(x_0)}{\partial x^L \partial x^M} \,
    \times \notag \\
    && \quad \times b^{L,j_2}(x_0) \, b^{M,j_2}(x_0) \cdot (t-t_0)^2 \notag \\
    && + \frac{1}{4} \sum_{j_1, j_2 =1}^m \,\, \sum_{J,K,L,M=1}^d
    \frac{\partial^2 f(x_0)}{\partial x^J \partial x^K}
    \frac{\partial b^{J,j_1}(x_0)}{\partial x^L} \,
    b^{L,j_2}(x_0) \, \times \notag \\
    && \quad \times \frac{\partial b^{K,j_1}(x_0)}{\partial x^M}
    \, b^{M,j_2}(x_0) \cdot (t-t_0)^2
    + \mathcal{R}_2(t,t_0) . \notag
\end{eqnarray}
\begin{Bem}
    We have to point out that Theorem~\ref{Ito-tree-expansion-exact-sol:Wm}
    as well as Theorem~\ref{St-tree-expansion-exact-sol:Wm} provide an expansion
    of functionals of the solution of the corresponding deterministic
    differential equation system, i.e.\ SDE~(\ref{Intro-Ito-St-SDE1-integralform-Wm})
    with $b \equiv 0$. Here, trees composed of the root $\gamma$ and
    deterministic nodes $\tau$ have to be considered only.
    Further, the presented Taylor formulas coincide with the ones
    used for the numerical analysis of deterministic Runge-Kutta
    methods (see, e.g.\ \cite{Butcher87}).
\end{Bem}
To illustrate the stochastic Taylor expansion, we consider the
stochastic differential equation
\begin{equation}
    d X_t = \alpha X_t \, dt + \beta X_t \, dW_t, \qquad
    X_{0}=x_0 \in \mathbb{R}
\end{equation}
in the case of $d=m=1$, i.e.\ for a one-dimensional Wiener
process, with some constants $\alpha, \beta \in \mathbb{R}$. If we
choose $f(x)=x$, the expectation of the solution can be calculated
as $E^{0,x_0}(X_t)= x_0 \exp(\alpha \, t)$. In the case of $m=1$,
the single index $j_1$ at the trees can be omitted since we have
$j_1 \equiv 1$. By the determination of the elementary
differentials appearing in~(\ref{Example-allg-expansion}) we yield
the following truncated expansion
\begin{equation}
    \begin{split}
        E^{0,x_0}(X_t) = & x_0 + \alpha \, x_0 \, (t-0) + \tfrac{1}{2} \alpha^2
        \, x_0 \, (t-0)^2 + \mathcal{R}_2(t,0)
    \end{split}
\end{equation}
which coincides with the order two Taylor polynomial of the exact
solution. \\ \\
The Taylor expansion can also be calculated with the aid of
hirarchical sets and the operators $L^0$ and $L^j$ introduced by
Platen and Wagner~(see \cite{KP99,PlWa82}). However, the rooted
tree theory allows the direct determination of each elementary
differential by the corresponding tree which makes it more
transparent and allows an easy analysis of the structure of the
elementary differentials. For example, the rooted tree approach
serves as an appropriate tool for the determination of stochastic
Runge-Kutta methods (see \cite{BuBu00a,Roe03}).

\end{document}